\documentclass[11pt,onecoulme]{article}
\usepackage{amsfonts}
\usepackage{mathrsfs}
\usepackage{amssymb}
\usepackage{color, amsmath,amssymb, amsfonts, amstext,amsthm, latexsym}

\usepackage{amssymb, epsfig, amssymb, latexsym}
\usepackage{amsmath}
\usepackage{graphicx}
\usepackage{longtable}

\numberwithin{equation}{section}

\allowdisplaybreaks

\newtheorem{definition}{Definition}[section]
\newtheorem{theorem}[definition]{Theorem}
\newtheorem{lemma}[definition]{Lemma}
\newtheorem{remark}[definition]{Remark}

\newtheorem{hyp}[definition]{Hypothesis}

\title{Viscosity Solutions to First Order Path-Dependent Hamilton-Jacobi-Bellman Equations \thanks{This work was partially supported by  the National Natural Science Foundation of China  (Grant No. 11401474),  the  Natural Science Foundation of Shaanxi Province
               (Grant No. 2017JM1016)  and the Fundamental Research Funds for the Central Universities (Grant No. 2452019075).}}

\author{  Jianjun Zhou   \\
       College of Science,
             Northwest A\&F University,\\ Yangling 712100, Shaanxi, P. R.
             China\\
      \emph{E-mail:zhoujianjun@nwsuaf.edu.cn} }
      \date{}

\begin{document}

\maketitle

\pagestyle{plain}

\begin{abstract}
In this article, a notion of viscosity solutions is introduced for first order path-dependent Hamilton-Jacobi-Bellman (PHJB) equations associated with optimal control problems for path-dependent differential equations. We identify the value functional of the optimal control problems as the unique viscosity solution to the associated PHJB equations. We also show that our notion of viscosity solutions is consistent with the corresponding notion of classical solutions, and satisfies a stability property.

\medskip

\end{abstract}


\section{Introduction}
\par
            In the early 1980's,  Crandall and Lions \cite{cra1}
            introduced the notion of viscosity solutions to  first order Hamilton-Jacobi-Bellman (HJB) equations. Lions \cite{lio} applied this notion to deterministic optimal control problems. From then on,  a large number of  papers have been published
           developing
            the theory of viscosity solutions.  We refer to the survey paper
of Crandall, Ishii and Lions \cite{cran2}. Soon afterwards,  Crandall and Lions \cite{cra3}, \cite{cra4}, \cite{cra5}, \cite{cra6} and \cite{cra7}
  systematically introduced  the corresponding theory  for  viscosity
              solutions in
                  infinite dimensional Hilbert space. Then, Li and Yong \cite{li} studied the  general unbounded first-order HJB equations in infinite dimensional Hilbert space.
\par
 For the path-dependent  case, the theory of viscosity solutions is more difficult.
 Lukoyanov \cite{luk} developed
a theory of viscosity solutions to fully non-linear path-dependent first order  Hamilton-Jacobi equations.
The existence and uniqueness theorems
are proved when Hamilton function $\mathbf{H}$ is $d_p$-locally Lipschitz continuous in the path function.
For the second order path-dependent case, a viscosity solution
approach has been successfully initiated by Ekren, Keller, Touzi
 and Zhang \cite{ekren1} in the semilinear context, and further extended to
 fully nonlinear
equations by Ekren, Touzi, and Zhang \cite{ekren3, ekren4},  elliptic
equations by Ren \cite{ren}, obstacle problems by Ekren \cite{ekren0}, and degenerate second-order
equations by Ren, Touzi, and Zhang \cite{ren1}  and Ren and Rosestolato \cite{ren2}.
Cosso, Federico, Gozzi, Rosestolato, and Touzi \cite{cosso} studied a class of semilinear second order Path-dependent partial differential  equations with a linear unbounded operator on Hilbert
space.
%
%
\par
In this paper, we consider the following   controlled path-dependent   differential
                 equation:
\begin{eqnarray}\label{state1}
\begin{cases}
            \dot{X}^{\gamma_t,u}(s)=
           F(X_s^{\gamma_t,u},u(s)),  \ \ s\in [t,T],\\
~~~~X_t^{\gamma_t,u}=\gamma_t\in {\Lambda}_t.
\end{cases}
\end{eqnarray}
In the equation above,  the control $u(\cdot)$ belongs to
$
                             {\cal{U}}[t,T]:=\{u(\cdot):[t,T]\rightarrow U|\ u(\cdot) \ \mbox{is
                             measurable}\}
$
 with  $(U,d)$ is a   metric space;
 $T>0$ is an arbitrarily fixed finite time horizon;
 the unknown $X^{\gamma_t,u}$ is an $R^d$-valued function on $[0,T]$, denote by $X^{\gamma_t,u}(s)$ the value of $X^{\gamma_t,u}$  at
 time $s$, and $X^{\gamma_t,u}_s$ the whole history path of $X^{\gamma_t,u}$ from time 0 to $s$; $\Lambda_t$ denotes the set of all
 continuous $R^d$-valued functions defined over $[0,t]$ and ${\Lambda}=\bigcup_{t\in[0,T]}{\Lambda}_{t}$; $\gamma_t$ is an element of $\Lambda_t$ and denote by $\gamma_t(s)$ the value of $\gamma_t$  at
 time $s$.
We define a   norm on ${\Lambda}_t$  and a metric on ${\Lambda}$ as follows: for any $0\leq t\leq \bar{t}\leq T$ and $\gamma_t,\bar{\gamma}_{\bar{t}}\in {\Lambda}$,
\begin{eqnarray*}
   ||\gamma_t||_0:=\sup_{0\leq s\leq t}|\gamma_t(s)|,\ \ d_{\infty}(\gamma_t,\bar{\gamma}_{\bar{t}})
   :=|t-\bar{t}|
               +\sup_{0\leq s\leq \bar{t}}|\gamma_{t}(s\wedge t)-\bar{\gamma}_{\bar{t}}(s)|.
\end{eqnarray*}
 We assume the coefficient $F:\Lambda\times U\rightarrow R^d$   satisfies  Lipschitz condition under $||\cdot||_0$
                         with respect to  the path function.
\par
              We wish to minimize a cost functional of the form:
\begin{eqnarray}\label{cost1}
                     J(\gamma_t,u(\cdot)):=\int_{t}^{T}q(X^{\gamma_t,u}_\sigma,u(\sigma))d\sigma
                          +\phi(X^{\gamma_t,u}_T),\ \ \ (t,\gamma_t)\in [0,T]\times {\Lambda},
\end{eqnarray}
 over  ${{\mathcal
                  {U}}}[t,T]$.
Here $
        q: {\Lambda}\times U\rightarrow R$ and $\phi: {\Lambda}_T\rightarrow R
$   satisfy Lipschitz conditions under  $||\cdot||_0$ 
                         with respect to  the path function. 
             We define the value functional of the  optimal
                  control problem as follows:
\begin{eqnarray}\label{value1}
V(\gamma_t):=\inf_{u(\cdot)\in{\mathcal
                  {U}}[t,T]}J({\gamma_t,u(\cdot)}),\ \ \ (t,\gamma_t)\in [0,T]\times {\Lambda}.
\end{eqnarray}
  The goal of this article is to characterize this value functional  $V$. We
                         consider the following path-dependent
                     HJB (PHJB) equation:
  \begin{eqnarray}\label{hjb1}
\begin{cases}
\partial_tV(\gamma_t)+{\mathbf{H}}(\gamma_t,\partial_xV(\gamma_t))= 0,\ \ \  (t,\gamma_t)\in [0,T)\times {\Lambda},\\
 V(\gamma_T)=\phi(\gamma_T), \ \ \ \gamma_T\in {\Lambda}_T,
 \end{cases}
\end{eqnarray}
      where
$$
                                {\mathbf{H}}(\gamma_t,p)=\inf_{u\in{
                                         {U}}}[
                        (p,F(\gamma_t,u))_{R^d}  +q(\gamma_t,u)], \ \ \ (t,\gamma_t,p)\in [0,T]\times{\Lambda}\times R^d.
$$
Here    $(\cdot,\cdot)_{R^d}$ denotes the scalar product of $R^d$,
                         $ \partial_t$ and $\partial_x$ denote formally the so-called pathwise  (or functional or Dupire; see \cite{dupire1, cotn0, cotn1}) derivatives, where $ \partial_t$  is known as horizontal derivative, while $\partial_x$ is vertical derivative.
\par
                         The primary objective of this article is to
                         develop the concept  of  viscosity solutions
                         to
                          equation  (\ref{hjb1}). To focus on the main
idea, we address the Lipschitz case under $||\cdot||_0$.  
                          We shall show that the value functional
                         $V$  defined in  (\ref{value1}) is the   unique viscosity solution to
                         equation given in  (\ref{hjb1}).
\par
The main difficulty for our case
lies in both facts that the path space $\Lambda_T$ is an infinite dimensional Banach space, and that the maximal norm $||\cdot||_0$ is not G\^ateaux differentiable.
In order to study the PHJB equations defined in path space $\Lambda$, we need to give a suitable definition to ensure
that the value functional is a viscosity solution of the PHJB equations.
It is more important to guarantee the
uniqueness of the solutions. With our assumptions of  coefficients $F$, $q$ and $\phi$, the value functional is only Lipschitz continuous under $||\cdot||_0$  with respect to the path function, then the auxiliary
function in the proof of uniqueness should include the term
$||\cdot||_0^2$ or a functional which is equivalent to $||\cdot||_0^2$.
 The lack
of G${\hat{a}}$teaux differentiability of $||\cdot||_0^2$ makes the definition of viscosity solutions
more complex.
\par
 As mentioned above, the notion of viscosity solutions for second order PHJB equations has been developed by many authors. However, none of the results we know are directly applicable to our situation.
 In the  papers, Ekren,
Touzi and Zhang  
\cite{ekren3} and \cite{ekren4}, Ekren \cite{ekren0} and Ren \cite{ren}, in particular, the nondegeneracy assumption required on the  Hamilton function $\mathbf{H}$ is not fulfilled if the diffusion term is identically equal zero. In Ren,  Touzi and  Zhang \cite{ren1} and Ren and Rosestolato \cite{ren2}, the degenerate case is taken into consideration, but in order to apply these results one has to require that   the coefficients $F$, $q$ and $\phi$ are  $d_p$-uniformly continuous with respect to  the path function. We  notice that the first order PHJB equations were also studied in  Section 8 of Ekren,
Touzi and Zhang   \cite{ekren3}, and  the comparison principle was  established  when  the Hamilton function  $\mathbf{H}$ is  locally uniformly continuous under $d_\infty$ in  the path function. 
However, this condition is not satisfied in our case, as
our Hamilton function  ${\mathbf{H}}(\gamma_t,p)$ includes  the term $(p,F(\gamma_t,u))_{R^d}$ for every $(t,\gamma_t,p,u)\in [0,T]\times{\Lambda}\times R^d\times U$.

\par
Our main contributions are as follows. We want to extend the results in \cite{luk} to Lipschitz continuous case under $||\cdot||_0$.
This extension is nontrivial  since the maximal norm $||\cdot||_0$ is not G\^ateaux differentiable. 
     To overcome this difficulty,  we define a functional $\Upsilon^2:\Lambda\rightarrow R$ by
     $$\Upsilon^2(\gamma_t)=S(\gamma_t)+2|\gamma_t(t)|^2,\ \ (t,\gamma_t)\in [0,T]\times\Lambda,$$
     where
 \begin{eqnarray*}
S(\gamma_t)=\begin{cases}
            \frac{(||\gamma_{t}||_0^2-|\gamma_{t}(t)|^2)^2}{||\gamma_{t}||^2_0}, \
         ~~ ||\gamma_{t}||_0\neq0, \\
0, \ ~~~~~~~~~~~~~~~~~~~ ||\gamma_{t}||_0=0,
\end{cases} \ \ \ \  (t,\gamma_t)\in [0,T]\times\Lambda.
\end{eqnarray*}
This functional is the key to prove the uniqueness of viscosity solutions.
   We first show that it is equivalent to $||\cdot||_0^2$ and study its  regularity in the horizontal/vertical sense mentioned above. 
  We next obtain that $\Upsilon^2$ satisfies  a functional formula. This is important as  functional $\Upsilon^2$ is equivalent to $||\cdot||_0^2$. Then we can define an auxiliary function 
which includes  the  functional $\Upsilon^2$
 and  prove  the uniqueness of  viscosity
                 solutions.
  Regarding existence, 
we notice that the solution $X^{\gamma_t,u}(\cdot)$ to equation
        (\ref{state1}) is Lipshitz continuous with respect to the time $s\in [t,T]$, then as in Lukoyanov \cite{luk},
we can give a definition of viscosity solutions in a sequence of
bounded and uniformly Lipschitz continuous paths spaces
which are compact subsets of $\Lambda$, and prove that the value functional $V$ is a solution under our definition by 
dynamic programming principle.
\par
             The remaining of this article is organized as follows. In the following
              section,  we introduce some basic
             notations to be used throughout this paper, and prove
          Theorem \ref{theoremito} and  Lemmas \ref{theoremS},  \ref{theoremS000} and \ref{theoremsg0809827} which are
    the key of  the existence and uniqueness results of viscosity solutions.
    In Section 3, we introduce  preliminary results on path-dependent  optimal control problems. We  give  the dynamic programming principle,  which will be used in the following sections.
             In Section 4, we introduce our notion of viscosity solutions to
           equation (\ref{hjb1}) and  prove  that the value functional $V$ defined by (\ref{value1}) is a viscosity solution. We also show
             the consistency with the notion of classical solutions and the stability result.
               The uniqueness of viscosity solutions for equation (\ref{hjb1}) is proven in Section 5. Finally,   in the Appendix we present some additional results. 

\section{Preliminary work}  \label{RDS}
\par
Let 
 $T>0$ be a fixed number.  For each  $t\in[0,T]$,
          define
         $\hat{\Lambda}_t:=D([0,t];R^d)$ as  the set of c$\grave{a}$dl$\grave{a}$g  $R^d$-valued
         functions on $[0,t]$.
        We denote $\hat{\Lambda}^t=\bigcup_{s\in[t,T]}\hat{\Lambda}_{s}$  and let  $\hat{\Lambda}$ denote $\hat{\Lambda}^0$.
\par
A very important remark on the notations: as in  Dupire \cite{dupire1}, we will denote elements of $\hat{\Lambda}$ by lower
case letters and often the final time of its domain will be subscripted, e.g. $\gamma\in \hat{\Lambda}_t\subset \hat{\Lambda}$ will be
denoted by $\gamma_t$. Note that, for any $\gamma\in  \hat{\Lambda}$, there exists only one $t$ such that $\gamma\in  \hat{\Lambda}_t$. For any $0\leq s\leq t$, the value of
$\gamma_t$ at time $s$ will be denoted by $\gamma_t(s)$. Moreover, if
a path $\gamma_t$ is fixed, the path $\gamma_t|_{[0,s]}$, for $0\leq s\leq t$, will denote the restriction of the path  $\gamma_t$ to the interval
$[0,s]$.  We also point out that the space $\hat{\Lambda}^t$  does not possess an algebraic structure since $\gamma_s+\eta_l$ is not well defined  for each $\gamma_s,\eta_l\in \hat{\Lambda}^t$ when $s\neq l$.
        \par
        Following Dupire \cite{dupire1},  for $ x\in R^d,\gamma_t\in \hat{\Lambda}_t$, $0\leq t\leq \bar{t}\leq T$, we define $\gamma^x_{t}\in\hat{\Lambda}_t$ and $\gamma_{t,\bar{t}}\in \hat{\Lambda}_{\bar{t}}$ as
\begin{eqnarray*}
  \gamma^x_{t}(s)&=&\gamma_t(s),\ \ s\in [0,t); \ \ \ \gamma^x_{t}(t)=\gamma_t(t)+x;\\
  \gamma_{t,\bar{t}}(s)&=&\gamma_t(s),\ \ s\in [0,t];\ \ \ \gamma_{t,\bar{t}}(s)=\gamma_t(t), \ \ s\in (t,\bar{t}].
\end{eqnarray*}
 We define a   norm on  $\hat{\Lambda}_t$  and a metric on $\hat{\Lambda}$ as follows: for any $0\leq t\leq \bar{t}\leq T$ and $\gamma_t,\bar{\gamma}_{\bar{t}}\in \hat{\Lambda}$,
\begin{eqnarray}\label{2.1}
   ||\gamma_t||_0:=\sup_{0\leq s\leq t}|\gamma_t(s)|,\ \ d_{\infty}(\gamma_t,\bar{\gamma}_{\bar{t}})
   :=|t-\bar{t}|
               +\sup_{0\leq s\leq \bar{t}}|\gamma_{t,\bar{t}}(s)-\bar{\gamma}_{\bar{t}}(s)|.
\end{eqnarray}
It is clear that  $(\hat{\Lambda}_t, ||\cdot||_0)$ is a Banach space for every $t\in [0,T]$. Moreover,  from Lemma \ref{lemma2.1111} in the Appendix, it follows that
$(\hat{\Lambda}^t, d_{\infty})$ is a complete metric space. We also clarify  that, for every $t\in(0,T]$, the Banach space $(\hat{\Lambda}_t, ||\cdot||_0)$ is not
separable . This  does not create problems here, as the Banach space $(\hat{\Lambda}_t, ||\cdot||_0)$ is only used to study the regularity of functional $S$ which will be  defined in (\ref{sre}), and our optimal control problems and the associated PHJB equations are considered in the continuous path space.
 \par Now we define the pathwise  derivatives of Dupire \cite{dupire1}.
 \begin{definition}\label{definitionc0} (Pathwise derivatives)
       Let $t\in [0,T)$ and  $f:\hat{\Lambda}^t\rightarrow R$.
\begin{description}
        \item{(i)}  Given $(s,\gamma_s)\in [t,T)\times \hat{\Lambda}^t$, the horizontal derivative of $f$ at $\gamma_s$ (if the corresponding limit exists and is finite) is defined as
        \begin{eqnarray}\label{2.3}
               \partial_tf(\gamma_s):=\lim_{h\rightarrow0,h>0}\frac{1}{h}\left[f(\gamma_{s,s+h})-f(\gamma_s)\right].
\end{eqnarray}
For the final time $T$, the horizontal derivative of $f$ at $\gamma_T\in\hat{\Lambda}^t$ (if the corresponding limit exists and is finite) is defined as
$$
\partial_tf(\gamma_T):=\lim_{s<T, s\uparrow T}\partial_tf(\gamma_T|_{[0,s]}).
$$
If the above limit exists and is finite for every $(s,\gamma_s)\in [t,T]\times \Lambda^t$, the functional $\partial_tf:\hat{\Lambda}^t\rightarrow R$ is called the horizontal derivative of $f$ with domain $\hat{\Lambda}^t$.
       \par
      \item{(ii)}  Given $(s,\gamma_s)\in [t,T]\times \hat{\Lambda}^t$, the vertical derivative  of $f$ at $\gamma_s$ (if all the corresponding limits exist and are finite) is defined as
      \begin{eqnarray}\label{2.20jia}
 \partial_{x}f(\gamma_s):=(\partial_{x_1}f(\gamma_s),\partial_{x_2}f(\gamma_s),\ldots, \partial_{x_d}f(\gamma_s)),
\end{eqnarray}
where
        \begin{eqnarray}\label{2.2}
 \partial_{x_i}f(\gamma_s):=\lim_{h\rightarrow0}\frac{1}{h}\left[f(\gamma_s^{he_i})-f(\gamma_s)\right],\ \
 i=1,2,\ldots,d,
\end{eqnarray}
with $e_1,e_2,\ldots,e_d$ is the standard  orthonormal basis of $R^d$.
If all the above limits exist and are finite for every $(s,\gamma_s)\in [t,T]\times \Lambda^t$, the map $\partial_xf:=(\partial_{x_1}f,\partial_{x_2}f,\ldots, \partial_{x_d}f):\hat{\Lambda}^t\rightarrow R^d$ is called the vertical  derivative of $f$ with domain $\hat{\Lambda}^t$.
\par
We take the convention that $\gamma_s$ is column vector, but $\partial_{x}f$ denotes row vector.
\end{description}
\end{definition}
\begin{definition}\label{definitionc}
       Let  $t\in [0,T)$ and  $f:\hat{\Lambda}^t\rightarrow R$  be given.
\begin{description}
        \item{(i)}
                 We say $f\in C^0(\hat{\Lambda}^t)$ if $f$ is continuous in $\gamma_s$ on $\hat{\Lambda}^t$ under $d_{\infty}$. 
       \par
      \item{(ii)} We say $f\in C^{1}(\hat{\Lambda}^t)\subset C^{0}(\hat{\Lambda}^t)$ if  $\partial_{x_1}f$, $\partial_{x_2}f$, $\ldots$, $\partial_{x_d}f$ and $\partial_tf$ exist on $ \hat{\Lambda}^t$ and are in $C^0(\hat{\Lambda}^t)$.
\end{description}
\end{definition}
\par
 For each  $t\in[0,T]$, let $\Lambda_t:= C([0,t],R^d)$ be the set of all continuous $R^d$-valued functions defined over $[0,t]$. We denote ${\Lambda}^t=\bigcup_{s\in[t,T]}{\Lambda}_{s}$  and let  ${\Lambda}$ denote ${\Lambda}^0$.
 Clearly, $\Lambda:=\bigcup_{s\in[0,T]}{\Lambda}_{s}\subset\hat{\Lambda}$, and each $\gamma\in \Lambda$ can also be viewed as an element of $\hat{\Lambda}$. $(\Lambda_t, ||\cdot||_0)$ is a
 Banach space, and $(\Lambda^t,d_{\infty})$ is a complete metric space.
 $f:\Lambda^t\rightarrow R$ and $\hat{f}:\hat{\Lambda}^t\rightarrow R$ are called consistent
  on $\Lambda^t$ if $f$ is the restriction of $\hat{f}$ on $\Lambda^t$.
 For every $t\in [0,T]$,  $\mu>0$ and $M_0>0$, we also define ${\cal{C}}^{\mu}_{t,M_0}$ by
$$
                         {\cal{C}}^{\mu}_{t,M_0}:=\bigg{\{}\gamma_s\in {\Lambda}^t:  ||\gamma_s||_0\leq M_0,   \sup_{0\leq l<r\leq s}\frac{|\gamma_s(l)-\gamma_s(r)|}{|l-r|}\leq \mu(1+M_0)
                            \bigg{\}}.
$$
For simplicity, we let   ${\cal{C}}^{\mu}_{M_0}$ denote ${\cal{C}}^{\mu}_{t,M_0}$  when $t=0$.
\par
\par
  We remark that, following Dupire  \cite{dupire1}, we study PHJB equation (\ref{hjb1}) in the metric space $(\Lambda,d_\infty)$ in the present paper, while some literatures (for example, \cite{cosso} and \cite{ekren3}) study PHJB equations in a complete pseudometric space. The reason to do this is that it is convenient to define ${\cal{C}}^{\mu}_{t,M_0}$ in our framework, which will be used to define our notion of viscosity solutions.
\begin{definition}\label{definitionc2}
       Let  $t\in [0,T)$ and  $f:{\Lambda}^t\rightarrow R$  be given.
\begin{description}
        \item{(i)}
                 We say $f\in C^0({\Lambda}^t)$ if $f$ is continuous in $\gamma_s$  on $\Lambda^t$ under $d_{\infty}$. 
\par
       \item{(ii)}  We say $f\in C^{1}({\Lambda}^t)\subset C^0({\Lambda}^t)$ if
                     there exists $\hat{f}\in C^{1}(\hat{{\Lambda}}^t)$ which is consistent with $f$ on $\Lambda^t$. 
\end{description}
\end{definition}

\par
%
               The following
                theorem  is needed to prove the existence  of viscosity solutions.
\begin{theorem}\label{theoremito}
\ \
Suppose $X$ is a  continuous function on $[0,T]$ and an absolutely continuous function on $[\hat{t}, T]$, and $u\in C^{1}({\Lambda}^{\hat{t}})$ for some $\hat{t}\in [0,T)$. Then for any $t\in [\hat{t},T]$:
\begin{eqnarray}\label{statesop0}
                 f(X_t)=f(X_{\hat{t}})+\int_{\hat{t}}^{t}\partial_tf(X_s)ds+\int^{t}_{\hat{t}}\partial_xf(X_s)dX(s),\ \ t\in [\hat{t},T].
\end{eqnarray}
 Here and in the following, for every $s\in [0,T]$, $X(s)$ denotes  the value of $X$  at
 time $s$, and $X_s$ the whole history path of $X$ from time 0 to $s$.
\end{theorem}
The proof is  similar to Theorem 4.1 in  Cont \& Fournie \cite{cotn1} (see also Dupire \cite{dupire1}).
For the convenience of readers, here we give its proof.
\par
{\bf  Proof}. \ \ For any  $t\in[\hat{t},T]$, denote $X^n(s)=X(s){\mathbf{1}}_{[0,\hat{t})}(s)+\sum^{2^n-1}_{i=0}X(t_{i+1}){\mathbf{1}}_{[t_i,t_{i+1})}+X(t){\mathbf{1}}_{\{t\}}(s)$, $s\in [0,t]$. 
Here $t_i=\hat{t}+\frac{i(t-\hat{t})}{2^n}$. For every $(s,\gamma_s)\in [0,T]\times\hat{\Lambda}$, define $\gamma_{s-}\in \hat{\Lambda}$ by
$$
               \gamma_{s-}(\theta)=\gamma_{s}(\theta),\ \ \theta\in [0,s),\ \ \mbox{and}\  \ \gamma_{s-}(s)=\lim_{\theta\uparrow s}\gamma_{s}(\theta).
$$ We start with the decomposition
\begin{eqnarray}\label{decom}
                   f({X^n_{t_{i+1}}}_{-})-f({X^n_{t_{i}}}_{-})=f({X^n_{t_{i+1}}}_{-})-f({X^n_{t_{i}}})
                   +f({X^n_{t_{i}}})-f({X^n_{t_{i}}}_{-}).
\end{eqnarray}
 Let $\psi(l)=f({X^n_{t_{i},t_i+l}})$, we have $f({X^n_{t_{i+1}}}_{-})-f({X^n_{t_{i}}})=\psi(h)-\psi(0)$, where $h=\frac{t-\hat{t}}{2^n}$. Let $\psi_{t^+}$ denote the right derivative of $\psi$, then
 $$
 \psi_{t^+}(l)=\lim_{\delta>0,\delta\rightarrow0}\frac{\psi(l+\delta)-\psi(l)}{\delta}=\lim_{\delta>0,\delta\rightarrow0}\frac{f({X^n_{t_{i},t_i+l+\delta}})-f({X^n_{t_{i},t_i+l}})}{\delta}
 =\partial_tf(X_{t_i,t_i+l}^n),\ \ l\in[0,h].
 $$
 By
 $$
 d_\infty({X^n_{t_{i},t_i+l_1}},{X^n_{t_{i},t_i+l_2}})=|l_1-l_2|, \ \ \ l_1, l_2\in [0,h],
$$
 and  $f\in C^{1}(\hat{\Lambda}^{\hat{t}})$, we have $\psi$ and $\psi_{t^+}$ is continuous on $[0,h]$, therefore,
 $$
                                f({X^n_{t_{i+1}}}_{-})-f({X^n_{t_{i}}})=\psi(h)-\psi(0)=\int^{h}_{0}\psi_{t^+}(l)dl=\int^{t_{i+1}}_{t_i}\partial_tf(X_{t_i,l}^n)dl, \ i\geq0.
$$
The term $f({X^n_{t_{i}}})-f({X^n_{t_{i}}}_{-})$ in (\ref{decom}) can be written $\pi(X(t_{i+1})-X(t_i))-\pi(0)$, where
$\pi(x)=f({X^n_{t_{i}}}_{-}+x{\mathbf{1}}_{\{t_{i}\}})$. Since $f\in C^{1}(\hat{\Lambda}^{\hat{t}})$, $\pi$ is a $C^1$ function and $\pi'(x)=\partial_xf({X^n_{t_{i}}}_{-}+x{\mathbf{1}}_{\{t_{i}\}})$.
Thus, we have that:
$$
                              \pi({X(t_{i+1})-X(t_i)})-\pi(0)=\int^{t_{i+1}}_{t_i}\partial_xf({X^n_{t_{i}}}_{-}+(X(s)-X(t_i)){\mathbf{1}}_{\{t_{i}\}})dX(s), \ i\geq 1.
$$
 Summing over $i\geq 0$ and denoting $i(s)$ the index such that $s\in [t_{i(s)},t_{i(s)+1})$, we obtain
$$
                                 f(X^n_t)-f(X^n_{\hat{t}})
                                =\int^{t}_{\hat{t}}\partial_tf(X_{t_{i(s),s}}^n)ds
                                 +\int^{t\vee t_1}_{t_1}\partial_xf({X^n_{t_{i(s)}}}_{-}+(X(s)-X(t_{i(s)})){\mathbf{1}}_{\{t_{i(s)}\}})dX(s).
$$
  $f(X^n_t)$ 
  converges to $f(X_t)$.
   Since all approximations of $X$ appearing in the various integrals have a $||\cdot||_{0}$-distance from $X_s$ less than $||X^n_s-X_s||_0\rightarrow0$, $f\in C^{1}(\hat{\Lambda}^{\hat{t}})$ implies that the integrands appearing in the above integrals converge respectively to $\partial_tf(X_s)$ and  $\partial_xf(X_s)$ as $n\rightarrow\infty$. By $X$ is continuous and $f\in C^{1}(\hat{\Lambda}^{\hat{t}})$, the integrands in the various above integrals are bounded.  The dominated convergence then ensure that the Lebesgue-Stieltjes integrals converge to the terms appearing in (\ref{statesop0}) as $n\rightarrow\infty$.\ \ $\Box$
  \par
By the above Theorem, we have the following important results.
\begin{lemma}\label{0815lemma}
   For every $t\in[0,T)$,
             let $f\in C^{1}(\Lambda^t)$ and $\hat{f}\in C^{1}(\hat{\Lambda}^t)$ such that $\hat{f}$ is consistent with $f$ on $\Lambda^t$, then the following definition
             $$
             \partial_tf:=\partial_t\hat{f}, \ \ \ \partial_xf:=\partial_x\hat{f}, \ \ \mbox{on} \ \Lambda^t
             $$
              is independent of the choice of $\hat{f}$. Namely, if there is another $\hat{f}'\in C^{1}(\hat{\Lambda}^t)$ such that $\hat{f}'$ is consistent with $f$ on $\Lambda^t$, then the derivatives of $\hat{f}'$
              coincide with those of $\hat{f}$ on $\Lambda^t$.
\end{lemma}
\par
{\bf  Proof}. \ \ By the definition of the horizontal derivative, it is clear that $\partial_t\hat{f}(\gamma_s)=\partial_t\hat{f}'(\gamma_s)$ for every $(s,\gamma_s)\in [t,T]\times\Lambda^t$.
We claim that $\partial_x\hat{f}(\gamma_s)=\partial_x\hat{f}(\gamma_s)$ also holds  for every $(s,\gamma_s)\in [t,T]\times\Lambda^t$. In fact, if not, there exist a constant $\varepsilon>0$ and $(s,\gamma_s)\in [t,T]\times \Lambda^t$ such that
$$
     |\partial_x\hat{f}(\gamma_s)-\partial_x\hat{f}(\gamma_s)|\geq \varepsilon>0.
$$
 Without loss of generality, we may assume
 $$
     |\partial_{x_1}\hat{f}(\gamma_s)-\partial_{x_1}\hat{f}(\gamma_s)|\geq \varepsilon>0.
$$
 By the  regularity $\partial_{x_1}\hat{f},\partial_{x_1}\hat{f}'\in C^0(\hat{\Lambda}^t)$,  we can assume $s<T$. Define $X:[0,T]\rightarrow R^d$ by
 $$
       X(\sigma)=\gamma_s(\sigma), \ \ \sigma\in [0,s], \ \ \ X(\sigma)=\gamma_s(s)+(\sigma-s)e_1, \ \ \sigma\in (s,T].
 $$
 It is clear that $X$ is  continuous  on $[0,T]$ and  absolutely continuous on $[s, T]$. Then by  Theorem \ref{theoremito},
\begin{eqnarray}\label{0823}
                        \int^{l}_{s}\partial_{x_1}\hat{f}(X_\sigma)d\sigma=\int^{l}_{s}\partial_{x_1}\hat{f}'(X_\sigma)d\sigma, \ \ l\in [s,T].
\end{eqnarray}
 On the other hand,  by also the  regularity $\partial_{x_1}\hat{f},\partial_{x_1}\hat{f}'\in C^0(\hat{\Lambda}^t)$, there exists a constant $\hat{l}\in (s,T]$ such that
 $$
  |\partial_{x_1}\hat{f}(X_\sigma)-\partial_{x_1}\hat{f}'(X_\sigma)|\geq \frac{\varepsilon}{2},\ \ \mbox{for all}\ \ \sigma\in [s,\hat{l}].
 $$
 Therefore,
 $$
                        \int^{l}_{s}|\partial_{x_1}\hat{f}(X_\sigma)-\partial_{x_1}\hat{f}'(X_\sigma)|d\sigma\geq \frac{\varepsilon}{2}|l-s|,   \ \ l\in [s,\hat{l}],
$$
 which contradict to  (\ref{0823}).
 \ \ $\Box$
                          \par
We conclude this section with  the following four  lemmas which will be used to prove the  stability and uniqueness  of viscosity solutions.
\begin{lemma}\label{lemmaleft}
(see Proposition 1  in  \cite{tang1})  For $t\in[0,T]$, $\mu>0$  and $M_0>0$, ${\cal{C}}^{\mu}_{t,M_0}$ is a compact subset of $(\Lambda^t,d_{\infty})$.
\end{lemma}
%
\par
For every fixed $(\hat{t},a_{\hat{t}})\in [0,T] \times \hat{\Lambda}_{\hat{t}}$, define $g^{a_{\hat{t}}}:\hat{\Lambda}^{\hat{t}}\rightarrow R$ by
$$g^{a_{\hat{t}}}(\gamma_t):=||\gamma_{t}-a_{\hat{t}, t}||_H^2,\ \ (t, \gamma_t)\in [\hat{t},T]\times\hat{\Lambda}^{\hat{t}},$$
where
$$
||{\gamma}_{{t}}||_H^2=\int^{t}_{0}|{\gamma}_{{t}}(s)|^2ds,\ \ \ (t,\gamma_t)\in[0,T]\times\hat{\Lambda}.
$$
\par
 We also define $S:\hat{\Lambda}\times \hat{\Lambda}\rightarrow R$ by 
 \begin{eqnarray}\label{sre}
S(\gamma_t,\gamma'_t)=\begin{cases}
            \frac{(||\gamma_{t}-\gamma'_{t}||_0^2-|\gamma_{t}(t)-\gamma'_{t}(t)|^2)^2}{||\gamma_{t}-\gamma'_{t}||^2_0}, \
         ~~ ||\gamma_{t}-\gamma'_{t}||_0\neq0; \\
0, \ ~~~~~~~~~~~~~~~~~~~~~~~~~~~~~~ ||\gamma_{t}-\gamma'_{t}||_0=0,
\end{cases} \ \ \  (t,\gamma_t), (t,\gamma'_t)\in [0,T]\times \hat{\Lambda}.
\end{eqnarray}
For simplicity, we let $S(\gamma_t)$ denote $S(\gamma_t,\gamma'_t)$ when $\gamma'_t(l)\equiv{0}$ for all $l\in [0,t]$.
\begin{lemma}\label{lemmagp2} For every fixed $(\hat{t},a_{\hat{t}})\in [0,T) \times \hat{\Lambda}_{\hat{t}}$,
 $g^{a_{\hat{t}}}\in C^{1}(\hat{\Lambda}^{\hat{t}})$.
\end{lemma}
   {\bf  Proof}. \ \
    It is clear that $g^{a_{\hat{t}}}\in C^0(\hat{\Lambda}^{\hat{t}})$ and $\partial_xg^{a_{\hat{t}}}(\gamma_s)=0$ for all $(s,\gamma_s)\in [\hat{t},T]\times\hat{\Lambda}^{\hat{t}}$.
   Now we consider $\partial_tg^{a_{\hat{t}}}$. For every $
                            (s, \gamma_s)\in [{\hat{t}},T)\times\hat{\Lambda}^{\hat{t}}$,
     \begin{eqnarray*}
   \partial_tg^{a_{\hat{t}}}(\gamma_s)=\lim_{h\rightarrow0,h>0}\frac{g^{a_{\hat{t}}}(\gamma_{s,s+h})-g^{a_{\hat{t}}}(\gamma_s)}{h}
    =\lim_{h\rightarrow0,h>0}\frac{\int^{s+h}_{s}|\gamma_{s}(s)-a_{\hat{t}}(\hat{t})|^2d\sigma}{h}=|\gamma_{s}(s)-a_{\hat{t}}(\hat{t})|^2.
   \end{eqnarray*}
   For $\gamma_T\in \hat{\Lambda}_T$,
   $$
\partial_tg^{a_{\hat{t}}}(\gamma_T):=\lim_{s<T,s\uparrow T}\partial_tg^{a_{\hat{t}}}(\gamma_T|_{[0,s]})=\left|\lim_{s\rightarrow T}\gamma_T(s)-a_{\hat{t}}(\hat{t})\right|^2.
$$
It is clear that $\partial_tg^{a_{\hat{t}}}\in C^0(\hat{\Lambda}^{\hat{t}})$. Thus, we show that $g^{a_{\hat{t}}}\in C^{1}(\hat{\Lambda}^{\hat{t}})$. \ \ $\Box$ 
\begin{lemma}\label{theoremS}
For every fixed $(\hat{t},a_{\hat{t}})\in [0,T) \times \hat{\Lambda}_{\hat{t}}$, define $S^{a_{\hat{t}}}:\hat{\Lambda}^{\hat{t}}\rightarrow R$ by
 $$S^{a_{\hat{t}}}(\gamma_t):=S(\gamma_{t},a_{\hat{t}, t}),\ \ (t,\gamma_t)\in [\hat{t}, T]\times\hat{\Lambda}^{\hat{t}}.$$
 Then $S^{a_{\hat{t}}}\in C^{1}(\hat{\Lambda}^{\hat{t}})$. Moreover,
\begin{eqnarray}\label{s0}
                \frac{3-5^{\frac{1}{2}}}{2}||\gamma_t||_0^2\leq   S(\gamma_t)+   |\gamma_t(t)|^2\leq 2||\gamma_t||_0^2, \ (t,\gamma_t)\in [0,T]\times\hat{\Lambda}.
\end{eqnarray}
\end{lemma}
\par
   {\bf  Proof  }. \ \  First, by the definition of $S^{a_{\hat{t}}}$, it is clear that $S^{a_{\hat{t}}}\in C^0(\Lambda^{\hat{t}})$ and $\partial_tS^{a_{\hat{t}}}(\gamma_{t})=0$ for all $(t,\gamma_t)\in [\hat{t},T]\times \hat{\Lambda}^{\hat{t}}$.
  Second, we consider $ \partial_{x_i}S^{a_{\hat{t}}}$. For every $(t,\gamma_t)\in [\hat{t},T]\times\hat{\Lambda}^{\hat{t}}$,
   \begin{eqnarray*}
   &&\partial_{x_i}S^{a_{\hat{t}}}(\gamma_{t})=\lim_{h\rightarrow0}\frac{S^{a_{\hat{t}}}(\gamma_t^{he_i})-S^{a_{\hat{t}}}(\gamma_{t})}{h}
   =\lim_{h\rightarrow0}\frac{S(\gamma_t^{he_i},{a_{\hat{t},t}})-S(\gamma_{t},{a_{\hat{t},t}})}{h}.
   \end{eqnarray*}
For every $(t,\gamma_t)\in [0,T]\times\hat{\Lambda}$, let $||\gamma_t||_{0^-}=\sup_{0\leq s<t}|\gamma_t(s)|$  and $\gamma^i_t(t)=\gamma_t(t)e_i,\ i=1,2,\ldots, d$.
Then,
 if $|\gamma_t(t)-a_{\hat{t}}(\hat{t})|<||\gamma_t-a_{\hat{t},t}||_{0^-}$,
\begin{eqnarray}\label{s1}
   \partial_{x_i}S^{a_{\hat{t}}}(\gamma_{t})
   &=&\lim_{h\rightarrow0}\frac{{(||\gamma_t-a_{\hat{t},t}||^2_{0}-|\gamma_t(t)+{he_i}-a_{\hat{t}}(\hat{t})|^2)^2}
   -{(||\gamma_t-a_{\hat{t},t}||^2_{0}-|\gamma_t(t)-a_{\hat{t}}(\hat{t})|^2)^2}}{h||\gamma_t-a_{\hat{t},t}||^2_{0}}\nonumber\\
   &=&-\frac{4(||\gamma_t-a_{\hat{t},t}||^2_{0}-|\gamma_t(t)-a_{\hat{t}}(\hat{t})|^2)(\gamma^i_t(t)-a^i_{\hat{t}}(\hat{t}))}{||\gamma_t-a_{\hat{t},t}||^2_{0}};
   \end{eqnarray}
 if  $|\gamma_t(t)-a_{\hat{t}}(\hat{t})|>||\gamma_t-a_{\hat{t},t}||_{0^-}$,
\begin{eqnarray}\label{s2}
   \partial_{x_i}S^{a_{\hat{t}}}(\gamma_{t})=0;
   \end{eqnarray}
if  $|\gamma_t(t)-a_{\hat{t}}(\hat{t})|=||\gamma_t-a_{\hat{t},t}||_{0^-}\neq0$,
since
\begin{eqnarray*}
&&||\gamma^{he_i}_t-a_{\hat{t},t}||_0^2-|\gamma_t(t)+{he_i}-a_{\hat{t}}(\hat{t})|^2\\
&=&
\begin{cases}
0,\ \ \ \ \ \ \ \ \ \  \ \ \ \ \ \ \ \ \ \ \  \ \ \ \ \ \ \ \ \ \  \ \ \ \ \ \ \ \ \ \ \ \ \ \ \ \ \ \ |\gamma_t(t)+he_i-a_{\hat{t}}(\hat{t})|\geq |\gamma_t(t)-a_{\hat{t}}(\hat{t})|,\\
  |\gamma_t(t)-a_{\hat{t}}(\hat{t})|^2-|\gamma_t(t)+{he_i}-a_{\hat{t}}(\hat{t})|^2,\ \ |\gamma_t(t)+he_i-a_{\hat{t}}(\hat{t})|<|\gamma_t(t)-a_{\hat{t}}(\hat{t})|,
  \end{cases}
\end{eqnarray*}
we have
\begin{eqnarray}\label{s3}
 0\leq\lim_{h\rightarrow0}\bigg{|}\frac{S^{a_{\hat{t}}}(\gamma_t^{he_i})-S^{a_{\hat{t}}}(\gamma_t)}{h}\bigg{|}
  \leq\lim_{h\rightarrow0}\bigg{|}\frac{h^2{(h+2(\gamma^i_t(t)-a^i_{\hat{t}}(\hat{t})))^2}}{h||\gamma^{he_i}_t-a_{\hat{t},t}||_0^2} \bigg{|}=0;
   \end{eqnarray}
    if $|\gamma_t(t)-a_{\hat{t}}(\hat{t})|=||\gamma_t-a_{\hat{t},t}||_{0^-}=0$,
\begin{eqnarray}\label{ss4}
  \partial_{x_i}S^{a_{\hat{t}}}(\gamma_{t})=0.
   \end{eqnarray}
From (\ref{s1}), (\ref{s2}), (\ref{s3}) and (\ref{ss4}) we obtain that, for all $(t,\gamma_t)\in  [\hat{t},T]\times\Lambda^{\hat{t}}$,
\begin{eqnarray}\label{0810zhou}
    \partial_{x_i}S^{a_{\hat{t}}}(\gamma_t)=\begin{cases}-\frac{4(||\gamma_t-a_{\hat{t},t}||^2_{0}-|\gamma_t(t)-a_{\hat{t}}(\hat{t})|^2)
    (\gamma^i_t(t)-a^i_{\hat{t}}(\hat{t}))}{||\gamma_t-a_{\hat{t},t}||^2_{0}}, \  \ \  ||\gamma_t-a_{\hat{t},t}||^2_{0}\neq0,\\
    0, ~~~~~~~~~~~~~~~~~~~~~~~~~~~~~~~~~~~~~~~~~~~~~~~~||\gamma_t-a_{\hat{t},t}||^2_{0}=0.
    \end{cases}
\end{eqnarray}
It is clear that $\partial_{x_i}S^{a_{\hat{t}}}\in C^0(\hat{\Lambda}^{\hat{t}})$.
 Thus, we have show that $S^{a_{\hat{t}}}\in C^{1}(\hat{\Lambda}^{\hat{t}})$.
\par
Now we prove (\ref{s0}). It is clear that
$$
                      S(\gamma_t)+   |\gamma_t(t)|^2\leq 2||\gamma_t||_0^2, \ (t, \gamma_t)\in [0,T]\times \hat{\Lambda}.
$$
On the other hand, for every $(t, \gamma_t)\in [0,T]\times \hat{\Lambda}$,
$$
                      S(\gamma_t)+   |\gamma_t(t)|^2\geq \frac{3-5^{\frac{1}{2}}}{2}||\gamma_t||_0^2, \ \  \mbox{if }\ \ ||\gamma_t||_0^2-|\gamma_t(t)|^2\leq\frac{5^{\frac{1}{2}}-1}{2}||\gamma_t||_0^2,
$$
 and 
$$
                      S(\gamma_t)+   |\gamma_t(t)|^2
                      \geq \frac{(\frac{5^{\frac{1}{2}}-1}{2})^2||\gamma_t||_0^4}{||\gamma_t||_0^2}=\frac{3-5^{\frac{1}{2}}}{2}||\gamma_t||_0^2\ \ \mbox{if }\ \ ||\gamma_t||_0^2-|\gamma_t(t)|^2>\frac{5^{\frac{1}{2}}-1}{2}||\gamma_t||_0^2.
$$
Thus, we  have (\ref{s0}) holds true.
 The proof is now complete. \ \ $\Box$
\par
For every constant $M>0$, define
$$
         \Upsilon^M(\gamma_t):=S(\gamma_t)+M|\gamma_t(t)|^2, \ \  (t, \gamma_t)\in [0,T]\times \hat{\Lambda};
$$
and
$$
         \Upsilon^M(\gamma_t,\eta_s)= \Upsilon^M(\eta_s,\gamma_t):=\Upsilon^M(\eta_s-\gamma_{t,s}), \ \  0\leq t\leq s\leq T, \ \gamma_t,\eta_s\in \hat{\Lambda}.
$$
 The following lemma will be used to prove Lemma \ref{lemma4.4}.
\begin{lemma}\label{theoremS000} For $M\geq2$, we have that
\begin{eqnarray}\label{jias5}
2\Upsilon^M(\gamma_t)+2\Upsilon^M(\gamma'_t)\geq\Upsilon^M(\gamma_t+\gamma'_t), \ \  (t, \gamma_t, \gamma'_t)\in [0,T]\times \hat{\Lambda}\times \hat{\Lambda}.
\end{eqnarray}
\end{lemma}
\par
   {\bf  Proof  }. \ \  If one of $||\gamma_t||_0$, $||\gamma'_t||_0$ and $||\gamma_t+\gamma'_t||_0$ is equal to $0$, it is clear that (\ref{jias5}) holds. Then we may assume that
      all of $||\gamma_t||_0$, $||\gamma'_t||_0$ and $||\gamma_t+\gamma'_t||_0$ are not equal to $0$.
      By the definition of $\Upsilon^M$, we get, for every $(t, \gamma_t, \gamma'_t)\in [0,T]\times \hat{\Lambda}\times \hat{\Lambda}$,
   \begin{eqnarray*}
   \Upsilon^M(\gamma_t+\gamma'_t)
   =||\gamma_t+\gamma'_t||_0^2+\frac{|\gamma_t(t)+\gamma'_t(t)|^4}{||\gamma_t+\gamma'_t||_0^2}+(M-2)|\gamma_t(t)+\gamma'_t(t)|^2.
\end{eqnarray*}
Letting $x:=||\gamma_t+\gamma'_t||_0^2$ and $a:=|\gamma_t(t)+\gamma'_t(t)|^2$, we have
$$
   \Upsilon^M(\gamma_t+\gamma'_t)=f(x,a):=x+\frac{a^2}{x}+(M-2)a.
$$
By $$
              f_x(x,a)=1-\bigg{(}\frac{a}{x}\bigg{)}^2\geq0,\ \  f_a(x,a)=2\frac{a}{x}+M-2\geq0,\ \ \ \forall \ x>0, \ a\geq0,
$$
and
$$||\gamma_t+\gamma'_t||_0^2\leq 2||\gamma_t||^2+2||\gamma'_t||_0^2,\ \ |\gamma_t(t)+\gamma'_t(t)|^2\leq 2|\gamma_t(t)|^2+2|\gamma'_t(t)|^2,
$$
 we obtain that
$$
   \frac{1}{2}\Upsilon^M(\gamma_t+\gamma'_t)\leq||\gamma_t||^2+||\gamma'_t||_0^2
   +\frac{(|\gamma_t(t)|^2+|\gamma'_t(t)|^2)^2}{||\gamma_t||^2+||\gamma'_t||_0^2}+(M-2)(|\gamma_t(t)|^2+|\gamma'_t(t)|^2).
$$
Combining with
$$
   \Upsilon^M(\gamma_t)+\Upsilon^M(\gamma'_t)=||\gamma_t||_0^2+||\gamma'_t||_0^2
   +\frac{|\gamma_t(t)|^4}{||\gamma_t||_0^2}+\frac{|\gamma'_t(t)|^4}{||\gamma'_t||_0^2}
   +(M-2)(|\gamma_t(t)|^2+|\gamma'_t(t)|^2),
$$
 we have
$$
   \Upsilon^M(\gamma_t)+\Upsilon^M(\gamma'_t)-\frac{1}{2}\Upsilon^M(\gamma_t+\gamma'_t)
   \geq\frac{|\gamma_t(t)|^4}{||\gamma_t||_0^2}+\frac{|\gamma'_t(t)|^4}{||\gamma'_t||_0^2}
   -\frac{(|\gamma_t(t)|^2+|\gamma'_t(t)|^2)^2}{||\gamma_t||^2+||\gamma'_t||_0^2}.
$$
Let $c=\frac{|\gamma_t(t)|^2}{||\gamma_t||_0^2}$, $b=\frac{|\gamma'_t(t)|^2}{||\gamma'_t||_0^2}$, $z=||\gamma_t||_0^2$ and $y=||\gamma'_t||_0^2$,
 we get that
\begin{eqnarray*}
   &&(||\gamma_t||_0^2+||\gamma'_t||_0^2)[\Upsilon^M(\gamma_t)+\Upsilon^M(\gamma'_t)-\frac{1}{2}\Upsilon^M(\gamma_t+\gamma'_t)]\\
   &\geq&(z+y)(c^2z+b^2y)-(cz+by)^2=(c-b)^2zy\geq0.
\end{eqnarray*}
Thus we obtain (\ref{jias5}) holds true.
The proof is now complete. \ \ $\Box$
\par
 Combing Theorem \ref{theoremito} and Lemmas \ref{0815lemma}, \ref{lemmagp2} and  \ref{theoremS}, we obtain
 \begin{lemma}\label{theoremsg0809827} For every fixed $(\hat{t},a_{\hat{t}})\in [0,T) \times {\Lambda}_{\hat{t}}$, if restrict $g^{a_{\hat{t}}}$ and $S^{a_{\hat{t}}}$ on ${\Lambda}^{\hat{t}}$ still denoted by themselves, then $g^{a_{\hat{t}}}\in C^{1}({\Lambda}^{\hat{t}})$ and $S^{a_{\hat{t}}}\in C^{1}({\Lambda}^{\hat{t}})$. Moreover, if $X$ is a  continuous function on $[0,T]$ and an absolutely continuous function on $[\hat{t}, T]$, we have
\begin{eqnarray*}
g^{a_{\hat{t}}}(X_s)=g^{a_{\hat{t}}}(X_{\hat{t}})+\int^{s}_{\hat{t}}|X(\sigma)-a_{\hat{t}}(\hat{t})|^2d\sigma, \  \ s\in [\hat{t},T];
\end{eqnarray*}
and
$$
S^{a_{\hat{t}}}(X_s)=S^{a_{\hat{t}}}(X_{\hat{t}})+\int^{s}_{\hat{t}}\partial_xS^{a_{\hat{t}}}(X_\sigma)dX(\sigma), \  \ s\in [\hat{t},T],
$$
where
$$
\partial_xS^{a_{\hat{t}}}=(\partial_{x_1}S^{a_{\hat{t}}},\partial_{x_2}S^{a_{\hat{t}}},\ldots, \partial_{x_d}S^{a_{\hat{t}}})
$$
with $\partial_{x_i}S^{a_{\hat{t}}}$ is defined in (\ref{0810zhou}) for all $ i=1,2,\ldots,d$.
\end{lemma}

\section{A DPP for optimal control problems.}
\par
In this section, we consider the controlled state
             equation (\ref{state1}) and  value functional  (\ref{value1}).
 Let  $(U,d)$ is a metric space. An admissible control  $u(\cdot):=\{u(r),  r\in [t,s]\}$ on $[t,s]$ (with $0\leq t\leq s\leq T$) is a  measurable function taking values in $U$. The set of all admissible controls on $[t,s]$ is denoted by ${\cal{U}}[t,s]$,
      i.e.,
$$
                             {\cal{U}}[t,s]:=\{u(\cdot):[t,s]\rightarrow U|\ u(\cdot) \ \mbox{is
                             measurable}\}.
$$
Let us make the
following assumptions.
\begin{hyp}\label{hypstate}
       $ F:{\Lambda}\times U\rightarrow R^d$, $
        q: {\Lambda}\times U\rightarrow R$ and $\phi: {\Lambda}_T\rightarrow R
$  are continuous, such that for some  constant $L>0$,  for all $(t,\gamma_t,u)$,  $ (t, \gamma'_{t},u) \in [0,T]\times {\Lambda}\times U$ and $\eta_T, \eta'_T \in  {\Lambda}_T$,
                  \begin{eqnarray*}
                |F(\gamma_t,u)-F(\gamma'_{t},u)|&\vee&|q(\gamma_t,u)-q(\gamma'_{t},u)| \leq
                 L||\gamma_t-\gamma'_{t}||_0;\\
                  |F(\gamma_t,u)|\vee|q(\gamma_t,u)| &\leq&
                 L(1+||\gamma_t||_0);\\
                  |\phi(\eta_T)-\phi(\eta'_T)|&\leq& L||\eta_T-\eta'_{T}||_0;\\
                |\phi(\eta_T)|&\leq& L(1+||\eta_T||_0).
\end{eqnarray*}
\end{hyp}
We remark that Hypothesis \ref{hypstate} dose not include the linear-quadratic case. Here, to focus on the main
idea, we
prefer to work with the Lipschitz  case  and leave the linear-quadratic case for
future study.
\par
The following theorem is standard, 
 but we do not find it in the existing literature. For the convenience of readers, we give its proof.
%
\par
\begin{theorem}\label{theorem2.2} \ \ Assume that Hypothesis \ref{hypstate}  holds. Then for every $u(\cdot)\in {\cal{U}}[t,T]$
and $(t,\gamma_t)\in [0,T]\times {\Lambda}$, equation  (\ref{state1}) admits a
unique  solution $X^{\gamma_t,u}$.
 Moreover,
\begin{eqnarray}\label{2.3}
                  \sup_{s\in [0,T]}|X^{\gamma_t,u}(s)|\leq C_1{(}1+||\gamma_t||_{0}{)},
\end{eqnarray}
              where the constant $C_1$  depends only on $L$ and  $T$.
\end{theorem}
{\bf  Proof  }. \ \
We define a mapping $\Phi$ from $C([0,T];R^d)$ to itself by the formula
\begin{eqnarray*}
                 \Phi(X)(s)=\gamma_t(t)+\int^{s}_{t}F(X_\sigma,u(\sigma))d\sigma,\  s\in [t,T], \ \ \  \Phi(X)(s)=\gamma_t(s),\ \ s\in [0,t),
\end{eqnarray*}
and show that it is a contraction, under an equivalent norm $||X||=\sup_{s\in[0,T]}e^{-\beta s}|X(s)|$, where $\beta>0$ will be chosen later. By Hypothesis \ref{hypstate}
\begin{eqnarray*}
                 ||\Phi(X)||&\leq& ||\gamma_t||_0+L\sup_{s\in[t,T]}\left[e^{-\beta s}\int^{s}_{t}(1+||X_\sigma||_0)d\sigma\right]\\
                 &\leq&  ||\gamma_t||_0+LT+L\sup_{s\in[t,T]}\int^{s}_{t}e^{-\beta (s-\sigma)}e^{-\beta \sigma}||X_\sigma||_0d\sigma\\
                 &\leq& ||\gamma_t||_0+LT+\frac{L}{\beta}||X||.
\end{eqnarray*}
 This show that $\Phi$ is a well defined mapping on $C([0,T];R^d)$. If $X,X'$ are functions belonging to this space, similar passages show that
 \begin{eqnarray*}
                 ||\Phi(X)-\Phi(X')||\leq \frac{L}{\beta}||X-X'||.
\end{eqnarray*}
Therefore, for $\beta>L$, the mapping is a contraction. In particular, we obtain $||X^{\gamma_t,u}||\leq C_1(1+||\gamma_t||_0)$, which prove the
estimate (\ref{2.3}).\ \ $\Box$
\par
               Let us now consider the continuous dependence of the solution $X^{\gamma_t,u}(\cdot)$ to equation
        (\ref{state1}) on the initial condition, the property will be used in the proof of Theorem \ref{theorem3.2222}.
\begin{theorem}\label{theorem2.3}\ \ Assume that Hypothesis \ref{hypstate}  holds. Then,  constant $C_2>0$ exists that depend only on $L$ and $T$,
   such that, for every  $0\leq t_1\leq t_2\leq T$, and  $\gamma_{t_1}^1,\gamma_{t_2}^2\in{\Lambda}$,
\begin{eqnarray}\label{2.4}
            \sup_{u(\cdot)\in {\cal{U}}[t_1,T]}\left|\left|X^{\gamma_{t_1}^1,u}_T-X^{\gamma_{t_2}^2,u}_T\right|\right|_0
              \leq C_2\left[\left|\left|\gamma_{t_1,t_2}^1-\gamma_{t_2}^2\right|\right|_0+\left(1+\left|\left|\gamma^1_{t_1}\right|\right|_0\right)(t_2-t_1)\right].
\end{eqnarray}
\end{theorem}
   \par
{\bf  Proof}. \ \
                    For any $0\leq t_1\leq t_2\leq T$ and $\gamma_{t_1}^1,\gamma_{t_2}^2\in{\Lambda}$, let $X^{u,i}_s$ denote $X^{\gamma_{t_i}^i,u}_s$ for $s\in [t_i,T]$, where $i=1,2$. Thus, we
                obtain
$$
                  \left|\left|X^{u,1}_l-X^{u,2}_l\right|\right|_0
                \leq
                \left|\left|\gamma_{t_1,t_2}^1-\gamma_{t_2}^2\right|\right|_0+L\left(1+\left|\left|X^{u,1}_{t_2}\right|\right|_0\right)(t_2-t_1)
               +L\int^{l}_{t_2}\left|\left|X^{u,1}_\sigma-X^{u,2}_\sigma\right|\right|_0d\sigma.
$$
               Using the Gronwall-Bellman inequality, by (\ref{2.3}),  we
                         obtain the following result, for a constant $C_2>0$  depending only on $L$ and $T$,
$$
                  \left|\left|X^{u,1}_T-X^{u,2}_T\right|\right|_0
                  \leq
                  C_2\left[\left|\left|\gamma_{t_1,t_2}^1-\gamma_{t_2}^2\right|\right|_0+\left(1+\left|\left|\gamma^1_{t_1}\right|\right|_0\right)(t_2-t_1)\right].
$$
              Applying the supremum i.e.,  $\sup_{u(\cdot)\in {\cal{U}}[t_1,T]}$, to both sides of the previous inequality, we get (\ref{2.4}).
         \ \ $\Box$
         \par
         The following  theorem  show that the solution $X^{\gamma_t,u}(\cdot)$ to equation
        (\ref{state1}) is Lipshitz continuous with respect to the time $s\in [t,T]$  even if  the initial value $(t,\gamma_t)$ belongs to $[0,T]\times\Lambda$. The result will be used to prove the existence of viscosity solutions in Theorem \ref{theoremvexist}.
\begin{theorem}\label{theorem2.33333333}\ \ Assume that Hypothesis \ref{hypstate}  holds. Then,  constant $C_3>0$ exists that depend only on $L$ and $T$,
   such that, for every  $(t,\gamma_{t})\in[0,T]\times{\Lambda}$,
\begin{eqnarray}\label{2.6}
            \sup_{u(\cdot)\in {\cal{U}}[t,T]}|X^{\gamma_t,u}(s_2)-X^{\gamma_t,u}(s_1)|\leq
            C_3(1+||\gamma_t||_0)|s_2-s_1|, \ \ \ t\leq s_1\leq s_2\leq T.
\end{eqnarray}
\end{theorem}
   \par
{\bf  Proof}. \ \
 For any $0\leq t\leq s_1\leq s_2\leq T$ and $\gamma_t\in\Lambda$, by (\ref{2.3}), we
                obtain the following result:
\begin{eqnarray*}
                   |X^{\gamma_t,u}(s_2)-X^{\gamma_t,u}(s_1)|
               \leq L(1+C_1(1+||\gamma_t||_0))|s_2-s_1|.
\end{eqnarray*}
    Taking the supremum in ${\cal{U}}[t,T]$, we obtain (\ref{2.6}).
         \ $\Box$
\par
                Our first result about the value functional   is the local boundedness  and
                two kinds of
                continuities.
                \begin{theorem}\label{theorem3.2222} Suppose that Hypothesis \ref{hypstate} holds true. Then, there exists a constant $C_4>0$  
such that,
                          for every $0\leq t\leq s\leq T$ and  $\gamma_t, \eta_t, \gamma'_s\in \Lambda$,
\begin{eqnarray}\label{3.5}
             |V(\gamma_t)|\leq C_4(1+||\gamma_t||_0);\ \            |V(\gamma_t)-V(\eta_t)|\leq
                        C_4
                        ||\gamma_t-\eta_t||_0;
\end{eqnarray}
\begin{eqnarray}\label{3.6}
                        &&|V(\gamma_t)-V(\gamma'_s)|\leq
                        C_4(1+||\gamma_t||_0\vee||\gamma'_s||_0)d_\infty{(}\gamma_t,\gamma'_s).
\end{eqnarray}
\end{theorem}
 \par
{\bf  Proof}. \ \
                  By Hypothesis \ref{hypstate}, (\ref{2.3}) and (\ref{2.4}), for any $u(\cdot)\in {\cal{U}}[t,T]$,
                  we have
                  \begin{eqnarray*}
                        &&|J(\gamma_t,u(\cdot))-J(\gamma'_{s},u(\cdot))|\\
                        &\leq& L\int_{t}^{s}(1+||X^{\gamma_t,u}_\sigma||_0)d\sigma+L(T+1)
                        ||X^{\gamma_t,u}_T-X^{\gamma'_{s},u}_T||_0\\
                        &\leq&
                        L(T+1)C_2
                        (||\gamma_{t,s}-\gamma'_s||_0+
                        (1+||\gamma_t||_0)(s-t))
                        +L(1+C_1(1+||\gamma_t||_0))(s-t).
\end{eqnarray*}
                  Thus, taking the infimum in $u(\cdot)\in {\cal{U}}[t,T]$, we can find a constant $C_4>0$
                  such that  (\ref{3.6}) holds.
                   By the similar procedure, we can show
                     (\ref{3.5}) holds true. The theorem is proved. \ \ $\Box$
\par
                        Now we present the dynamic programming
                        principle (DPP) for optimal control problems (\ref{state1}) and (\ref{value1}).
\begin{theorem}\label{theorem3.3}
                      Assume the Hypothesis \ref{hypstate}  holds true. Then, for every $(t,\gamma_t)\in [0,T)\times \Lambda$ and $s\in [t,T]$, we have that
\begin{eqnarray}\label{3.7}
             V(\gamma_t)=\inf_{u(\cdot)\in {\cal{U}}[t,T]}\bigg{[}\int_{t}^{s}q(X^{\gamma_t,u}_\sigma,u(\sigma))d\sigma+V(X^{\gamma_t,u}_s)\bigg{]}.
\end{eqnarray}
\end{theorem}
The proof is very similar to the case without path-dependent (see Theorem 2.1 in page 160 of \cite {yong}). For the convenience of readers, here we give its proof.
 \par
{\bf  Proof}. \ \ First of all, for any  $u(\cdot)\in{\cal{U}}[s,T]$,
$s\in[t,T]$ and any  $u(\cdot)\in{\cal{U}}[t,s]$, by putting them
             concatenatively, we get $u(\cdot)\in{\cal{U}}[t,T]$.
                  Let us denote the right-hand side of (\ref{3.7}) by
                  $\overline{V}({\gamma_t})$. By (\ref{value1}), we have
$$
           V({\gamma_t})\leq J({\gamma_t,u(\cdot)})
              =\int_{t}^{s}q(X^{\gamma_t,u}_\sigma,u(\sigma))d\sigma+J(X^{{\gamma_t,u}}_s,{{u(\cdot)}}),\
              u(\cdot)\in {\cal{U}}[t,T].
$$
             Thus, taking the infumum over $u(\cdot)\in
             {\cal{U}}[s,T]$, we obtain
$$
           V({\gamma_t})\leq \int_{t}^{s}q(X^{\gamma_t,u}_\sigma,u(\sigma))d\sigma+V(X^{{\gamma_t,u}}_s).
$$
             Consequently,
$$
                 V({\gamma_t})\leq \overline{V}({\gamma_t}).
$$
            On the other hand, for any $\varepsilon>0$, there exists a $u^\varepsilon(\cdot)\in {\cal{U}}[t,T]$ such that
       \begin{eqnarray*}
           V({\gamma_t})+\varepsilon&\geq& J({\gamma_t},u^\varepsilon(\cdot))
           =
              \int_{t}^{s}q(X^{\gamma_t,u^\varepsilon}_\sigma,{u^\varepsilon}(\sigma))d\sigma
                           +J(X^{\gamma_t,u^\varepsilon}_s,{u^\varepsilon}(\cdot))\\
              &\geq&\int_{t}^{s}q(X^{\gamma_t,u^\varepsilon}_\sigma,{u^\varepsilon}(\sigma))d\sigma
                           +V(X^{\gamma_t,u^\varepsilon}_s)\geq\overline{V}({\gamma_t}).
\end{eqnarray*}
             Hence, (\ref{3.7}) follows.\ \ $\Box$

\section{Viscosity solutions to  HJB equations: Existence theorem.}

\par
                    In this section, we consider the  first  order PHJB  equation (\ref{hjb1}). As usual, we start with classical solutions.
\par
\begin{definition}\label{definitionccc}     (Classical solution)
              A functional $v\in C^{1}({\Lambda})$       is called a classical solution to the PHJB equation (\ref{hjb1}) if it satisfies
                equation (\ref{hjb1}) point-wisely.
 \end{definition}
        We shall get that  the value functional $V$ defined by (\ref{value1}) is a viscosity solution of equation (\ref{hjb1}).  We  give the following definition for the viscosity solutions.
\par
  For every $M_0>0$, $\mu>0$, $(t,\gamma_t)\in [0,T]\times \Lambda$ and $w\in C^0(\Lambda)$, define
$$
             J^+_{\mu, M_0}(\gamma_t,w):=\bigg{\{}\varphi\in C^1({\Lambda}^t):  0=({w}-{{\varphi}})({\gamma_t})=\sup_{\eta_s\in {\cal{C}}^\mu_{t,M_0}}
                         ({w}- {{\varphi}})(
                         \eta_s)\bigg{\}},
$$
and
$$
             J^-_{\mu, M_0}(\gamma_t,w):=\bigg{\{}\varphi\in C^1({\Lambda}^t):  0=({w}+{{\varphi}})({\gamma_t})=\inf_{\eta_s\in {\cal{C}}^\mu_{t,M_0}}
                         ({w}+{{\varphi}})(
                         \eta_s)\bigg{\}}.
$$


\begin{definition}\label{definition4.1} \ \
  Let $w\in C^0({\Lambda})$.
\begin{description}
        \item{(i)}
         For any $\mu>0$,  $w$ is called  a
                            viscosity  $\mu$-subsolution (resp., $\mu$-supersolution)
                             of equation (\ref{hjb1}) if the terminal condition  $w(\gamma_T)\leq \phi(\gamma_T)$(resp., $w(\gamma_T)\geq \phi(\gamma_T)$),
                             $\gamma_T\in {\Lambda}_T$ is satisfied, and for every $M_0>0$,
                             whenever  $\varphi \in J^+_{\mu,M_0}(\gamma_s,w)$ (resp., $\varphi\in J^-_{\mu,M_0}(\gamma_s,w)$)
                                 with $(s,{\gamma}_{s})\in [0,T)\times {\cal{C}}^\mu_{M_0}$ and $|{\gamma}_{s}(s)|<
                   M_0$,  we have
$$
                          \partial_t{\varphi}({\gamma}_{s})
                           +{\mathbf{H}}({\gamma}_{s},\partial_x{\varphi}({\gamma}_{s}))\geq0,
$$
$$
                          (\mbox{resp.,}\ -\partial_t{\varphi}({\gamma}_{s})
                           +{\mathbf{H}}({\gamma}_{s},-\partial_x{\varphi}({\gamma}_{s})
                           )
                          \leq0).
$$
\par
       \item{(ii)}
        $w$ is called a
                            viscosity  subsolution (resp., supersolution)
                             of equation (\ref{hjb1}) if there exists a  $\mu_0>0$ such that, for all ${\mu}\geq \mu_0$, $w$ is a
                            viscosity  ${\mu}$-subsolution (resp., ${\mu}$-supersolution)
                             of equation (\ref{hjb1}).
\par
       \item{(iii)}
 $w\in C^0(\Lambda)$ is said to be a
                             viscosity solution of equation  (\ref{hjb1}) if it is
                             both a viscosity subsolution and a viscosity
                             supersolution.
\end{description}
%
%
\end{definition}
\begin{remark}\label{remarkv}
 Assume that  the coefficients $F(\gamma_t,u)=\overline{F}(t,\gamma_t(t),u)$,  $q(\gamma_t,u)=\overline{q}(t,\gamma_t(t),u)$ and
                      $\phi(\eta_T)=\overline{\phi}(\eta_T(T))$  for all
                     $(t,\gamma_t,u) \in [0,T]\times{\Lambda}\times U$ and $\eta_T\in \Lambda_T$. Then there exists a function $\overline{V}:[0,T]\times R^d\rightarrow R$ such that  $V(\gamma_t)=\overline{V}(t,\gamma_t(t))$ for all
                     $(t,\gamma_t) \in [0,T]\times{\Lambda}$, and the PHJB equation (\ref{hjb1}) reduces to the following HJB equation:
 \begin{eqnarray}\label{hjb3}
\begin{cases}
\overline{V}_{t^+}(t,x)+\overline{{\mathbf{H}}}(t,x,\nabla_x\overline{V}(t,x))= 0,\ \ \  (t, x)\in
                               [0,T)\times R^d,\\
 V(T,x)=\overline{\phi}(x), \ \ \ x\in R^d;
 \end{cases}
\end{eqnarray}
      where
$$
                                \overline{{\mathbf{H}}}(t,x,p)=\inf_{u\in{
                                         {U}}}[
                        (p,\overline{F}(t,x,u))_{R^d}  +\overline{q}(t,x,u)], \ \ \ (t,x,p)\in [0,T]\times R^d\times R^d.
$$
Here  and in the sequel, $\nabla_x$ denotes the standard first  order derivative
with respect to $x$. However, slightly different from the HJB literature, $\overline{V}_{t^+}$ denotes the right time-derivative of $\overline{V}$.
\end{remark}
\par
 The following theorem show that our definition of viscosity solutions to PHJB equation (\ref{hjb1}) is a natural  extension of classical viscosity solutions to HJB equation  (\ref{hjb3}).
 \begin{theorem}\label{theoremnatural} \ \ Consider the setting in Remark \ref{remarkv}.
                          Assume that $V$ is a viscosity solution of PHJB equation (\ref{hjb1}) in the sense of Definition \ref{definition4.1}. Then $\overline{V}$ is a
                          viscosity solution of HJB equation (\ref{hjb3}) in the standard sense (see Definition 2.4 on  page 165 of \cite{yong}).
\end{theorem}
 {\bf  Proof}. \ \ Without loss of generality, we shall only prove the viscosity subsolution property.
First,  from $V$ is a viscosity subsolution of equation (\ref{hjb1}), it follows that, for every $x\in R^d$,
 $$
                            \overline{V}(T,x)=V(\gamma_T)\leq \phi(\gamma_T)=\overline{\phi}(x),
$$
where $\gamma_T\in \Lambda$ with $\gamma_T(T)= x$.
\\
Next, let $\overline{\varphi}\in C^{1}([0,T]\times R^d)$ and $(t,x)\in [0,T)\times R^d$
 such that $$
                         0=(\overline{V}- {\overline{\varphi}})(t,x)=\sup_{(s,y)\in [0,T]\times R^d}
                         (\overline{V}- {\overline{\varphi}})(
                         s,y).
$$
 Define $\varphi:{\Lambda}\rightarrow R$ by
 $$
 \varphi(\gamma_s)=\overline{\varphi}(s,\gamma_s(s)),\  (s, \gamma_s)\in [0,T]\times\hat{\Lambda},
 $$
 and define $\hat{\gamma}_{t}\in \Lambda_t$ by
 $$
   \hat{\gamma}_{t}(s)= x,\ \ s\in [0,t].$$
 It is clear that $\varphi\in C^1({\Lambda})\subset C^1({\Lambda}^{t})$ and
 $$\partial_t\varphi(\gamma_s)=\overline{\varphi}_{t}(s,\gamma_s(s)),\ \ \partial_x\varphi(\gamma_s)=\nabla_x\overline{\varphi}(s,\gamma_s(s)),\ \ (s,\gamma_s)\in [0,T]\times{\Lambda}.
 $$
  Let $M_0>0$ be large enough such that
 $|x|<M_0$, since $\hat{\gamma}_{t}\in {\cal{C}}^\mu_{t,M_0}$  for all $\mu>0$,  by the definitions of $V$ and $\varphi$, we get that, for all $\mu>0$,
 $$
 0=(V-{{\varphi}})(\hat{\gamma}_{t})=(\overline{V}-\overline{\varphi})(t,x)=\sup_{(s,y)\in [0,T]\times R^d}
                         (\overline{V}- {\overline{\varphi}})(
                         s,y)=\sup_{\gamma_s\in {\cal{C}}^\mu_{t,M_0}}
                         (V-{{\varphi}})(
                         \gamma_s).  
 $$
 Therefore, for all $\mu>0$, we have  $ \varphi \in J^+_{\mu, M_0}(\hat{\gamma}_{t},V)$ with $(t,\hat{\gamma}_{t})\in [0,T)\times {\cal{C}}^\mu_{M_0}$ and $|\hat{\gamma}_{t}(t)|<M_0$.
 Since $V$ is a viscosity subsolution of PHJB equation (\ref{hjb1}), there exists a $\mu_0>0$ such that, for all ${\mu}\geq \mu_0$,
 $$\partial_t{\varphi}(\hat{\gamma}_{t})
                           +{\mathbf{H}}(\hat{\gamma}_{t},\partial_x{\varphi}(\hat{\gamma}_t))\geq0.
$$
Thus,
$$\overline{\varphi}_{t}(t,x)+\overline{{\mathbf{H}}}(t,x,\nabla_x\overline{\varphi}(t,x))\geq 0.$$
 By the arbitrariness of   $\overline{\varphi}\in C^{1}([0,T]\times R^d)$, we see that $\overline{V}$ is a viscosity subsolution of HJB equation (\ref{hjb3}), and thus completes the proof. \ \ $\Box$
\par
We are now in a  position  to give  the existence proof for  viscosity solutions.
\begin{theorem}\label{theoremvexist} \ \
                          Suppose that Hypothesis \ref{hypstate}  holds. Then the value
                          functional $V$ defined by (\ref{value1}) is a
                          viscosity solution to equation  (\ref{hjb1}).
\end{theorem}
 {\bf  Proof}. \ \         First, for every $M_0>0$,   let $\mu_0=C_3>0$. For every $\mu\geq \mu_0$,   we let  $\varphi\in J^+_{\mu,M_0}(\gamma_t,V)$
                  with
                   $(t,{\gamma}_t)\in [0,T)\times {\cal{C}}^\mu_{M_0}$ and $|{\gamma}_t(t)|<
                   M_0$. 
 For fixed $u\in U$,
by Theorem \ref{theorem2.33333333}, we can let $\delta>0$ be  small enough
 such that $t+\delta\leq T$, $||X_{t+\delta}^{\gamma_t,u}||_0\leq M_0$ and
 $$\sup_{t\leq s_1<s_2\leq t+\delta}\frac{|X^{\gamma_t,u}(s_2)-X^{\gamma_t,u}(s_1)|}{|s_2-s_1|}\leq
            C_3(1+||\gamma_t||_0)\leq \mu(1+M_0).
 $$
Combining with  ${\gamma}_t\in {\cal{C}}^\mu_{M_0}$, we have $X_{t+\delta}^{\gamma_t,u}\in {\cal{C}}^{\mu}_{t,M_0}$.
 Then by the DPP (Theorem \ref{theorem3.3}), we obtain that 
 \begin{eqnarray}\label{4.9}
                            {{\varphi}} (\gamma_t)=V(\gamma_t)
                           \leq \int^{t+\delta}_{t}
                              q(X_\sigma^{\gamma_t,u},u(\sigma))d\sigma
                              +V(X_{t+\delta}^{\gamma_t,u})
                              \leq \int^{t+\delta}_{t}
                              q(X_\sigma^{\gamma_t,u},u(\sigma))d\sigma
                              +{{\varphi}}(X_{t+\delta}^{\gamma_t,u}).
\end{eqnarray}
As $\varphi\in C^1(\Lambda^t)$, by Theorem \ref{theoremito} and Lemma \ref{0815lemma} we show that
  \begin{eqnarray*}
                 0  &\leq&\lim_{\delta\rightarrow0}\bigg{[}\frac{1}{\delta}\int^{t+\delta}_{t}
                              q(X_\sigma^{\gamma_t,u},u(\sigma))d\sigma
                            +\frac{1}{\delta}{[}{{\varphi}}(X_{t+\delta}^{\gamma_t,u})-{{\varphi}} (\gamma_t){]}\bigg{]}\\
                            &=&q(\gamma_t,u)+ \partial_t\varphi(\gamma_t)+(\partial_{x}\varphi(\gamma_t),
                     F(\gamma_t,u))_{R^d}.
\end{eqnarray*}
   Taking the minimum in $u\in U$, 
   we have
$$
0
    \leq\partial_t{\varphi}({\gamma}_{t})
                           +{\mathbf{H}}({\gamma}_{t},\partial_x{\varphi}({\gamma}_{t})).
$$
On the other hand, it is clear that $V(\gamma_T)\leq\phi(\gamma_T)$ for all $\gamma_T\in \Lambda_T$.
 Then  $V$ is a viscosity $\mu$-subsolution of equation (\ref{hjb1}) for all $\mu\geq \mu_0$. Thus  $V$ is a viscosity subsolution of equation (\ref{hjb1}).
\par
      Next,   for every $M_0>0$,   let $\mu_0=C_3>0$. For every $\mu\geq \mu_0$,   we let  $\varphi\in J^-_{\mu,M_0}(\gamma_t,V)$
                  with
                   $(t,{\gamma}_t)\in [0,T)\times {\cal{C}}^\mu_{M_0}$ and $|{\gamma}_t(t)|<
                   M_0$.
                   %
%
By Theorem \ref{theorem2.33333333}, we can  let $\delta>0$ be  small enough
 such that
 $t+\delta\leq T$, $\sup_{u(\cdot)\in {\cal{U}}[t,T]}||X_{t+\delta}^{\gamma_t,u}||_0\leq M_0$ and
 $$\sup_{u(\cdot)\in {\cal{U}}[t,T]}\sup_{t\leq s_1<s_2\leq t+\delta}\frac{|X^{\gamma_t,u}(s_2)-X^{\gamma_t,u}(s_1)|}{|s_2-s_1|}\leq
            C_3(1+||\gamma_t||_0)\leq \mu(1+M_0).
 $$
Combining with  ${\gamma}_t\in {\cal{C}}^\mu_{M_0}$, we have $X_{t+\delta}^{\gamma_t,u}\in {\cal{C}}^{\mu}_{t,M_0}$ for all $u(\cdot)\in {\cal{U}}[t,T]$. Then, for any $\varepsilon>0$,  by  the DPP (Theorem \ref{theorem3.3}), one can find a control  ${u}^\varepsilon(\cdot)\equiv u^{\varepsilon,\delta}(\cdot)\in {\cal{U}}[t,T]$ such
   that 
\begin{eqnarray*}
    \varepsilon\delta
    \geq \int^{t+\delta}_{t}
                              q(X_\sigma^{\gamma_t,{u}^\varepsilon},{u}^\varepsilon(\sigma))d\sigma
                              +V(X_{t+\delta}^{\gamma_t,{u}^\varepsilon})-V(\gamma_t)
                              \geq \int^{t+\delta}_{t}
                              q(X_\sigma^{\gamma_t,{u}^\varepsilon},{u}^\varepsilon(\sigma))d\sigma
                              -\varphi(X_{t+\delta}^{\gamma_t,{u}^\varepsilon})+\varphi(\gamma_t).
\end{eqnarray*}
                 Then, 
                  applying Theorem  \ref{theoremito}  and Lemma \ref{0815lemma} to $\varphi$,
                  we obtain that
\begin{eqnarray*}
                           \varepsilon&\geq& \frac{1}{\delta}\int^{t+\delta}_{t}
                              q(X_\sigma^{\gamma_t,{u}^\varepsilon},{u}^\varepsilon(\sigma))d\sigma
                                -\frac{\varphi(X_{t+\delta}^{\gamma_t,{u}^\varepsilon})-\varphi(\gamma_t)}{\delta}
                                 \\
           &=& -\partial_t\varphi(\gamma_t)+\frac{1}{\delta}\int^{t+\delta}_{t}[q(\gamma_t,{u}^\varepsilon(\sigma)) -(\partial_{x}\varphi(\gamma_t),
                     F(\gamma_t,{u}^\varepsilon(\sigma)))_{R^d} ]d\sigma+o(1)\\
            &\geq& -\partial_t\varphi(\gamma_t)
            +\inf_{u\in U}[q(\gamma_t,{u}) -(\partial_{x}\varphi(\gamma_t),
                     F(\gamma_t,{u}))_{R^d}]+o(1).
\end{eqnarray*}
Letting $\delta\downarrow 0$ and  $\varepsilon\rightarrow0$, we show that
$$
            0\geq
                      -\partial_t\varphi(\gamma_t) +{\mathbf{H}}({\gamma}_{t},-\partial_x{\varphi}({\gamma}_{t})).
$$
Moreover, we also have $V(\gamma_T)\geq\phi(\gamma_T)$ for all $\gamma_T\in \Lambda_T$.
    Therefore, $V$ is also a viscosity $\mu$-supsolution of (\ref{hjb1}) for all $\mu\geq \mu_0$. Thus  $V$ is a viscosity supsolution of equation (\ref{hjb1}).  This completes the proof.\ \ $\Box$
 \par
Now, let us give   the result of  classical  solutions, which show the consistency of viscosity solutions.
 \begin{theorem}\label{theorem3.2}
                      Let $V$ denote the value functional  defined by (\ref{value1}). If  $V\in C^{1}({\Lambda})$, then
                 $V$ is a classical solution of  equation (\ref{hjb1}).
\end{theorem}
{\bf  Proof }. \ \
First, using the  definition of $V$ yields  $V(\gamma_T)=\phi(\gamma_T)$ for all $\gamma_T\in \Lambda_T$. Next, for   fixed $(t,\gamma_t,u)\in [0,T)\times{\Lambda}\times U$,
                  from  the DPP (Theorem \ref{theorem3.3}), we obtain the following result:
 \begin{eqnarray}\label{4.9000}
                           && 0\leq \int_{t}^{t+\delta}q(X^{\gamma_t,u}_\sigma,u)d\sigma
                                +V(X^{\gamma_t,u}_{t+\delta})
                           -V(\gamma_t),\ \ 0<\delta< T-t.
\end{eqnarray}
               By  
               Theorem  \ref{theoremito}  and Lemma \ref{0815lemma},
                the  inequality above implies that
               \begin{eqnarray*}
    0&\leq&\lim_{\delta\rightarrow 0^+}\frac{1}{\delta}\bigg{[}\int_{t}^{t+\delta}q(X^{\gamma_t,u}_\sigma,u)d\sigma+V(X^{\gamma_t,u}_{t+\delta})-V(\gamma_t)\bigg{]}\\
                       &=&\partial_t V(\gamma_t)+ (
                        F(\gamma_t,u),\partial_xV(\gamma_t))_{R^d}+q(\gamma_t,u).
\end{eqnarray*}
     Taking the minimum in $u\in U$, we have that
\begin{eqnarray}\label{3.14}
    0\leq\partial_t V(\gamma_t)+ {\mathbf{H}}(\gamma_t,\partial_x V(\gamma_t)).
\end{eqnarray}
          On the other hand,  let $(t,\gamma_t)\in [0,T)\times\Lambda$ be fixed. Then, by   DPP (Theorem \ref{theorem3.3}) and  $V\in C^1(\Lambda)$,
           there exists an
                         $\tilde{u}(\cdot)\equiv u^{\varepsilon,\delta}(\cdot)\in {\cal{U}}[t,T]$ for any $\varepsilon>0$ and $0<\delta<T-t$ such that
\begin{eqnarray*}
                \varepsilon\delta &\geq& \int_{t}^{t+\delta}q(X^{{\gamma}_{{t}},{\tilde{u}}}_s,{\tilde{u}}(s))ds
                                +V(X^{{\gamma}_{{t}},{\tilde{u}}}_{{{t}}+\delta})
                 -V(\gamma_t)\\
                  &=& {\partial_t}V(\gamma_t)\delta+\int_{t}^{t+\delta}q(\gamma_t,\tilde{u}(\sigma))d\sigma+\bigg{(}{\partial_{x}V(\gamma_t)},\int^{t+\delta}_{t}
                        F(\gamma_t,{\tilde{u}}(\sigma))d\sigma\bigg{)}_{R^d}+o(\delta)\\
                     &\geq& {\partial_t}V(\gamma_t)\delta+{\mathbf{H}}(\gamma_t,\partial_x V(\gamma_t))\delta+o(\delta).
\end{eqnarray*}
Then, dividing through by $\delta$ and letting $\delta\rightarrow0^+$, we
obtain that
$$
                       \varepsilon\geq \partial_tV(\gamma_t)+{\mathbf{H}}(\gamma_t,\partial_{x}V(\gamma_t)).
$$
               The desired result  is obtained by combining the inequality given above with (\ref{3.14}). \ \ $\Box$
                                \par
We conclude this section with   the stability of viscosity solutions.
\begin{theorem}\label{theoremstability}
                      Let $\mu>0$, $F,q,\phi$ satisfy  Hypothesis \ref{hypstate}, and $v\in C^0(\Lambda)$. Assume
                       \item{(i)}      for any $\varepsilon>0$, there exist $F^\varepsilon,  q^\varepsilon, \phi^\varepsilon$ and $v^\varepsilon\in C^0(\Lambda)$ such that  $F^\varepsilon, q^\varepsilon, \phi^\varepsilon$ satisfy  Hypothesis \ref{hypstate} and $v^\varepsilon$ is a viscosity $\mu$-subsolution (resp., $\mu$-supsolution) of equation (\ref{hjb1}) with generators $F^\varepsilon,  q^\varepsilon, \phi^\varepsilon$;
                           \item{(ii)} as $\varepsilon\rightarrow0$, $(F^\varepsilon,  q^\varepsilon, \phi^\varepsilon,v^\varepsilon)$ converge to
                           $(F,  q, \phi, v)$ uniformly in the following sense: 
\begin{eqnarray}\label{sss}
                         \lim_{\varepsilon\rightarrow0}\sup_{(t,{\gamma_t},u)\in [0,T]\times\Lambda\times U}\sup_{\eta_T\in \Lambda_T}[(|F^\varepsilon-F|
                         +|q^\varepsilon-q|)({\gamma}_t,u)+|\phi^\varepsilon-\phi|(\eta_T)+|v^\varepsilon-v|({\gamma}_t)]=0.
\end{eqnarray}
                   Then $v$ is a viscosity $\mu$-subsoluiton (resp., $\mu$-supersolution) of equation (\ref{hjb1}) with generators $F,q,\phi$.
\end{theorem}
{\bf  Proof }. \ \ Without loss of generality, we shall only prove the viscosity subsolution property.
First,  from $v^{\varepsilon}$ is a viscosity $\mu$-subsolution of equation (\ref{hjb1}) with generators $F^{\varepsilon},  q^{\varepsilon}, \phi^{\varepsilon}$, it follows that
 $$
                            v^{\varepsilon}(\gamma_T)\leq \phi^{\varepsilon}(\gamma_T),\ \ \gamma_T\in \Lambda_T.
$$
Letting $\varepsilon\rightarrow0$, we  have
 $$
                            v(\gamma_T)\leq \phi(\gamma_T),\ \ \gamma_T\in \Lambda_T.
$$
Next, for every  $M_0>0$,  we let   $\varphi\in J^+_{\mu,M_0}(\hat{\gamma}_{\hat{t}}, v)$
                  with
  $(\hat{t},\hat{\gamma}_{\hat{t}})\in [0,T)\times{\cal{C}}^\mu_{M_0} $ and $|\hat{\gamma}_{\hat{t}}(\hat{t})|<M_0$.
 Denote $\varphi_{1}(\gamma_t):=\varphi(\gamma_t)+|t-\hat{t}|^2+||\gamma_{t}-\hat{\gamma}_{{\hat{t}}, t}||_H^2$ for all
 $(t,\gamma_t)\in [\hat{t},T]\times\Lambda^{\hat{t}}$. By Lemma \ref{theoremsg0809827}, we have  $\varphi_{1}\in C^{1}({\Lambda}^{\hat{t}})$.
 For every $\varepsilon>0$, from Lemma \ref{lemmaleft} it follows that there exists $({{t}_{\varepsilon}},{\gamma}^{\varepsilon}_{{t}_{\varepsilon}})\in [\hat{t},T]\times {\cal{C}}^\mu_{\hat{t},M_0}$ such that
$$
                         (v^{\varepsilon}- {{\varphi_{1}}})({{t}_{\varepsilon}},{\gamma}^{\varepsilon}_{{t}_{\varepsilon}})=\sup_{\gamma_s\in {\cal{C}}^\mu_{{\hat{t}}, M_0}}
                         (v^{\varepsilon}- {{\varphi_{1}}})(
                         \gamma_s).
$$
We claim that $d_{\infty}({\gamma}^{\varepsilon}_{{t}_{\varepsilon}},\hat{\gamma}_{\hat{t}})\rightarrow0$ as $\varepsilon\rightarrow0$.
 Indeed, if not, by  Lemma \ref{lemmaleft}, we may assume there exist  $({\bar{t}},\bar{\gamma}_{\bar{t}})\in [\hat{t},T]\times {\cal{C}}^\mu_{\hat{t},M_0}$  and a subsequence of $({{t}_{\varepsilon}},{\gamma}^{\varepsilon}_{{t}_{\varepsilon}})$ still denoted by themselves  such that $({\bar{t}},\bar{\gamma}_{\bar{t}})\neq(\hat{t},\hat{\gamma}_{\hat{t}})$ and $d_{\infty}({\gamma}^{\varepsilon}_{{t}_{\varepsilon}},\bar{\gamma}_{\bar{t}})\rightarrow0$
  as $\varepsilon\rightarrow0$. Thus
\begin{eqnarray*}
   &&(v- {{\varphi}})(\bar{\gamma}_{\bar{t}})=\lim_{\varepsilon\rightarrow0}(v- {{\varphi}})({\gamma}^{\varepsilon}_{{t}_{\varepsilon}})
   \leq(v- {{\varphi}})(\hat{\gamma}_{\hat{t}})=(v- {{\varphi}_1})(\hat{\gamma}_{\hat{t}})\\
   &=&\lim_{\varepsilon\rightarrow0}[(v- v^\varepsilon)(\hat{\gamma}_{\hat{t}})+ (v^{\varepsilon}- {{\varphi_{1}}})(\hat{\gamma}_{\hat{t}})]\leq \lim_{\varepsilon\rightarrow0}[(v- v^\varepsilon)(\hat{\gamma}_{\hat{t}})+ (v^{\varepsilon}- {{\varphi_{1}}})({\gamma}^{\varepsilon}_{{t}_{\varepsilon}})]\\
   &=&(v- {{\varphi}})(\bar{\gamma}_{\bar{t}})-|{\bar{t}}-\hat{t}|^2
   -||\gamma_{{\bar{t}}}-\hat{\gamma}_{{\hat{t}}, \bar{t}}||_H^2,
\end{eqnarray*}
 contradicting $|{\bar{t}}-\hat{t}|^2+||\gamma_{{\bar{t}}}-\hat{\gamma}_{{\hat{t}}, \bar{t}}||_H^2>0$.
Then, for any $\rho>0$, by (\ref{sss}) there exists $\varepsilon>0$ small enough such that
$$
            \hat{t}\leq {t}_{\varepsilon}< T,  \        \ |{\gamma}^{\varepsilon}_{{t}_{\varepsilon}}({t}_{\varepsilon})|< {M_0},\ \
             2|{t}_{\varepsilon}-\hat{t}|+|{\gamma}^{\varepsilon}_{{t}_{\varepsilon}}({t}_{\varepsilon})
                       -\hat{\gamma}_{{\hat{t}}}({\hat{t}})|^2\leq \frac{\rho}{4},$$
             and
             $$
|\partial_t{\varphi}({\gamma}^{\varepsilon}_{{t}_{\varepsilon}})-\partial_t{\varphi}(\hat{\gamma}_{\hat{t}})|\leq \frac{\rho}{4}, \ |I|\leq \frac{\rho}{4}, \ |II|\leq \frac{\rho}{4},
$$
where
$$
I={\mathbf{H}}^{\varepsilon}({\gamma}^{\varepsilon}_{{t}_{\varepsilon}},
                           \partial_x{\varphi}({\gamma}^{\varepsilon}_{{t}_{\varepsilon}}))
                           -{\mathbf{H}}({\gamma}^{\varepsilon}_{{t}_{\varepsilon}},
                           \partial_x{\varphi}({\gamma}^{\varepsilon}_{{t}_{\varepsilon}})),
$$
$$
II={\mathbf{H}}({\gamma}^{\varepsilon}_{{t}_{\varepsilon}},
                           \partial_x{\varphi}({\gamma}^{\varepsilon}_{{t}_{\varepsilon}}))
                       -{\mathbf{H}}(\hat{\gamma}_{\hat{t}},\partial_x{\varphi}(\hat{\gamma}_{\hat{t}})),
$$
and
$$
                                {\mathbf{H}}^{\varepsilon}(\gamma_t,p)=\inf_{u\in{
                                         {U}}}[
                        (p,F^{\varepsilon}(\gamma_t,u))_{R^d}
                        +q^{\varepsilon}(\gamma_t,u)],  \ \ (t,\gamma_t,p)\in [0,T]\times{\Lambda}\times R^d.
$$
 Since $v^{\varepsilon}$ is a viscosity $\mu$-subsolution of  equation (\ref{hjb1}) with generators $F^{\varepsilon}, q^{\varepsilon}, \phi^{\varepsilon}$, we have
$$
                          \partial_t{\varphi_{1}}({\gamma}^{\varepsilon}_{{t}_{\varepsilon}})
                           +{\mathbf{H}}^{\varepsilon}({\gamma}^{\varepsilon}_{{t}_{\varepsilon}},
                           \partial_x{\varphi_{1}}({\gamma}^{\varepsilon}_{{t}_{\varepsilon}}))\geq0.
$$
Notice that $\partial_x{\varphi_{1}}({\gamma}^{\varepsilon}_{{t}_{\varepsilon}})=\partial_x{\varphi}({\gamma}^{\varepsilon}_{{t}_{\varepsilon}})$, we obtain
\begin{eqnarray*}
                       0&\leq&   \partial_t{\varphi}({\gamma}^{\varepsilon}_{{t}_{\varepsilon}})
                       +2({t}_{\varepsilon}-\hat{t})+|{\gamma}^{\varepsilon}_{{t}_{\varepsilon}}({t}_{\varepsilon})
                       -\hat{\gamma}_{{\hat{t}}}({\hat{t}})|^2+{\mathbf{H}}(\hat{\gamma}_{\hat{t}},\partial_x{\varphi}(\hat{\gamma}_{\hat{t}}))\\
                       &&
                           +{\mathbf{H}}^{\varepsilon}({\gamma}^{\varepsilon}_{{t}_{\varepsilon}},
                           \partial_x{\varphi}({\gamma}^{\varepsilon}_{{t}_{\varepsilon}})) -{\mathbf{H}}({\gamma}^{\varepsilon}_{{t}_{\varepsilon}},
                           \partial_x{\varphi}({\gamma}^{\varepsilon}_{{t}_{\varepsilon}}))
                           +{\mathbf{H}}({\gamma}^{\varepsilon}_{{t}_{\varepsilon}},
                           \partial_x{\varphi}({\gamma}^{\varepsilon}_{{t}_{\varepsilon}}))
                           -{\mathbf{H}}(\hat{\gamma}_{\hat{t}},\partial_x{\varphi}(\hat{\gamma}_{\hat{t}}))\\
                       &\leq&\partial_t{\varphi}(\hat{\gamma}_{\hat{t}})+{\mathbf{H}}(\hat{\gamma}_{\hat{t}},
                       \partial_x{\varphi}(\hat{\gamma}_{\hat{t}}))+{\rho}.
\end{eqnarray*}
Letting $\rho\downarrow 0$, we show that
$$
\partial_t{\varphi}(\hat{\gamma}_{\hat{t}})+{\mathbf{H}}(\hat{\gamma}_{\hat{t}},
                       \partial_x{\varphi}(\hat{\gamma}_{\hat{t}}))\geq0.
$$
Since ${\varphi}\in C^{1}({\Lambda}^{{\hat{t}}})$ is arbitrary, we see that $v$ is a viscosity $\mu$-subsolution of equation (\ref{hjb1}) with generators $F,q,\phi$, and thus completes the proof.
\ \ $\Box$

\section{Viscosity solutions to  HJB equations: Uniqueness theorem.}
\par
             This section is devoted to a  proof of uniqueness of  viscosity
                   solutions to equation (\ref{hjb1}). This result, together with
                  the results from  the previous section, will be used to characterize
                   the value functional defined by (\ref{value1}).
                   \par
We  now state the main result of this section.
\begin{theorem}\label{theoremhjbm}  Suppose Hypothesis \ref{hypstate}   holds.
                         Let $W_1\in C^0({\Lambda})$ $(\mbox{resp}., W_2\in C^0({\Lambda}))$ be  a viscosity subsolution (resp., supsolution) to equation (\ref{hjb1}) and  let  there exist  constant $L>0$,
                        such that, for any  
                        $(t,\gamma_t),(s,\eta_s)\in[0,T]\times{\Lambda}$,
\begin{eqnarray}\label{w}
                                   |W_1(\gamma_t)|\vee |W_2(\gamma_t)|\leq L (1+||\gamma_t||_0);
                                   \end{eqnarray}
\begin{eqnarray}\label{w1}
                                |W_1(\gamma_t)-W_1(\eta_s)|\vee|W_2(\gamma_t)-W_2(\eta_s)|
                  \leq
                        L(1+||\gamma_t||_0\vee||\eta_s||_0)d_\infty{(}\gamma_t,\eta_s).
\end{eqnarray}
                   Then  $W_1\leq W_2$. 
\end{theorem}
\par
                      Theorems    \ref{theoremvexist} and \ref{theoremhjbm} lead to the result (given below) that the viscosity solution to the
                       PHJB equation given in (\ref{hjb1})
                      corresponds to the value functional  $V$ of our optimal control problem given in (\ref{state1}) and (\ref{value1}).
\begin{theorem}\label{theorem52}\ \
                 Assume that  Hypothesis \ref{hypstate}  holds. Then the value
                          functional $V$ defined by (\ref{value1}) is the unique viscosity
                          solution to equation (\ref{hjb1}) in the class of functionals satisfying (\ref{w}) and (\ref{w1}).
\end{theorem}
\par
   {\bf  Proof  }. \ \   Theorem \ref{theoremvexist} shows that $V$ is a viscosity solution to equation (\ref{hjb1}).
   Thus, our conclusion follows from  Theorems \ref{theorem3.2222} and
     \ref{theoremhjbm}.  \ \ $\Box$

\par
  Next, we prove Theorem \ref{theoremhjbm}.   Let $W_1$ be a viscosity subsolution of  equation (\ref{hjb1}). 
 We  note that for $\delta>0$, the functional
                    defined by $\tilde{W}:=W_1-\frac{\delta}{t}$ is a viscosity subsolution
                   for
 \begin{eqnarray*}
\begin{cases}
{\partial_t} \tilde{W}(\gamma_t)+{\mathbf{H}}(\gamma_t, \partial_x \tilde{W}(\gamma_t))
          = \frac{\delta}{t^2}, \ \  (t,\gamma_t)\in[0,T)\times {\Lambda}, \\
\tilde{W}(\gamma_T)=\phi(\gamma_T), \ \ \gamma_T\in \Lambda_T.
\end{cases}
\end{eqnarray*}
                As $W_1\leq W_2$ follows from $\tilde{W}\leq W_2$ in
                the limit $\delta\downarrow0$, it suffices to prove
                $W_1\leq W_2$ under the additional assumption given below:
$$
{\partial_t} {W_1}(\gamma_t)+{\mathbf{H}}(\gamma_t,  \partial_x {W_1}(\gamma_t))
          \geq c, \ \  (t,\gamma_t)\in[0,T)\times {\Lambda},\ \ c:=\frac{\delta}{T^2}.
$$
\par
   {\bf  Proof of Theorem \ref{theoremhjbm} } \ \   The proof of this theorem  is rather long. Thus, we split it into several
        steps.
\par
            $Step\  1.$ Definitions of auxiliary functions.
             \par
             By the definition of viscosity solutions, 
              there exists a $\mu_0>0$ such that $W_1$ $(\mbox{resp}., W_2)$ is a viscosity $\mu$-subsolution (resp.,  $\mu$-supsolution) to equation (\ref{hjb1}) for all $\mu\geq \mu_0$.
             \par
 We only need to prove that $W_1(\gamma_t)\leq W_2(\gamma_t)$ for all $(t,\gamma_t)\in
[T-\bar{a},T)\times
       {\Lambda}$.
        Here,
        $$\bar{a}=\frac{1}{96L}\wedge T.$$
         Then, we can  repeat the same procedure for the case
        $[T-i\bar{a},T-(i-1)\bar{a})$.  Thus, we assume the converse result that $(\breve{t},\breve{\gamma}_{\breve{t}})\in [T-\bar{a},T)\times
      {\Lambda}$ exists  such that
        $2\tilde{m}:=W_1(\breve{\gamma}_{\breve{t}})-W_2(\breve{\gamma}_{\breve{t}})>0$.  Because $\cup_{\mu\geq\mu_0,M>0}{\cal{C}}^\mu_{{\breve{t}},M}$  is dense in $\Lambda^{\breve{t}}$,
         by (\ref{w1}) there exist $\hat{\mu}\geq \mu_0, M_0>0$, $\tilde{t}\in [T-\bar{a},T)$ and $\tilde{\gamma}_{\tilde{t}}\in {\cal{C}}^{\hat{\mu}}_{\tilde{t},M_0}$
         such that
        $W_1(\tilde{\gamma}_{\tilde{t}})-W_2(\tilde{\gamma}_{\tilde{t}})>\tilde{m}$.
\par
       Let $
            \nu=1+\frac{1}{96TL},
$ and  consider that  $\varepsilon >0$ is  a small number such that
 $$
 W_1(\tilde{\gamma}_{\tilde{t}})-W_2(\tilde{\gamma}_{\tilde{t}})-2\varepsilon \frac{\nu T-\tilde{t}}{\nu
 T}(S(\tilde{\gamma}_{\tilde{t}})+|\tilde{\gamma}_{\tilde{t}}(\tilde{t})|^2)
 >\frac{\tilde{m}}{2},
 $$
      and
\begin{eqnarray}\label{5.3}
                          \frac{\varepsilon}{\nu T}\leq\frac{c}{2}.
\end{eqnarray}
  Next, we define for any $(t,\gamma_t,\eta_t)\in [T-\bar{a},T]\times{\Lambda}\times{\Lambda}$,
\begin{eqnarray*}
                 \Psi(\gamma_t,\eta_t)=W_1(\gamma_t)-W_2(\eta_t)-\frac{\alpha}{2}
                 \Upsilon^2(\gamma_t,\eta_t)-\varepsilon\frac{\nu T-t}{\nu
                 T}(S(\gamma_t)+S(\eta_t)+|\gamma_t(t)|^2+|\eta_t(t)|^2).
\end{eqnarray*}
 Finally, for the  fixed  ${\hat{\mu}}\geq\mu_0$  and every  $M\geq M_0$, 
we can apply Lemma \ref{lemmaleft} to find $({\hat{t}},\hat{{\gamma}}_{{\hat{t}}},\hat{{\eta}}_{{\hat{t}}})\in [T-\bar{a},T]\times
            {\cal{C}}_{M}^{\hat{\mu}}\times
            {\cal{C}}_{M}^{\hat{\mu}}$
          such that
 $$
            \Psi(\hat{{\gamma}}_{{\hat{t}}},\hat{{\eta}}_{{\hat{t}}})
                   \geq  \Psi(\tilde{\gamma}_{\tilde{t}},\tilde{\gamma}_{\tilde{t}})>\frac{\tilde{m}}{2} \ \ \mbox{and}   \ \       \Psi(\hat{{\gamma}}_{{\hat{t}}},\hat{{\eta}}_{{\hat{t}}})
                   \geq  \Psi(\gamma_t,\eta_t), \
                   (t,\gamma_t,\eta_t)\in [T-\bar{a},T]\times{\cal{C}}_{M}^{\hat{\mu}}\times{\cal{C}}_{M}^{\hat{\mu}}.
 $$
             We should note that the point
             $({\hat{t}},\hat{{\gamma}}_{{\hat{t}}},\hat{{\eta}}_{{\hat{t}}})$ depends on $\alpha,
             \hat{\mu}, \varepsilon,  M$.
 \par
 $Step\ 2.$ For the fixed  $\hat{\mu}\geq\mu_0$  and  every $M\geq M_0$, the following result  holds true:
 \begin{eqnarray}\label{5.10}
                       {\alpha}\Upsilon^2(\hat{{\gamma}}_{{\hat{t}}},\hat{{\eta}}_{{\hat{t}}})
                         \leq|W_1(\hat{{\gamma}}_{{\hat{t}}})-W_1(\hat{{\eta}}_{{\hat{t}}})|
                                   +|W_2(\hat{{\gamma}}_{{\hat{t}}})-W_2(\hat{{\eta}}_{{\hat{t}}})|\rightarrow0 \ 
                                   \mbox{as} \ \alpha\rightarrow+\infty.
 \end{eqnarray}
  Let us show the above. By the definition of $(\hat{{\gamma}}_{{\hat{t}}},\hat{{\eta}}_{{\hat{t}}})$, we have
 \begin{eqnarray}\label{5.56789}
                        2\Psi(\hat{{\gamma}}_{{\hat{t}}},\hat{{\eta}}_{{\hat{t}}})\geq  \Psi(\hat{{\gamma}}_{{\hat{t}}},\hat{{\gamma}}_{{\hat{t}}})
                        +\Psi(\hat{{\eta}}_{{\hat{t}}},\hat{{\eta}}_{{\hat{t}}}).
 \end{eqnarray}
 This implies that
 \begin{eqnarray}\label{5.6}
                         {\alpha}\Upsilon^2(\gamma_t,\eta_t)
                         &\leq&|W_1(\hat{{\gamma}}_{{\hat{t}}})-W_1(\hat{{\eta}}_{{\hat{t}}})|
                                   +|W_2(\hat{{\gamma}}_{{\hat{t}}})-W_2(\hat{{\eta}}_{{\hat{t}}})|\nonumber\\
                                   &\leq& 2L(2+||\hat{{\gamma}}_{{\hat{t}}}||_0+||\hat{{\eta}}_{{\hat{t}}}||_0)
                                   \leq 4L(1+M).
 \end{eqnarray}
   Letting $\alpha\rightarrow+\infty$, we get
                $$\Upsilon^2(\gamma_t,\eta_t)
                \rightarrow0\ 
                          \mbox{as} \ \alpha\rightarrow+\infty.
                          $$
Then from (\ref{s0}) it follows that
\begin{eqnarray}\label{5.66666123}
||\hat{{\gamma}}_{{\hat{t}}}-\hat{{\eta}}_{{\hat{t}}}||_0 \rightarrow0\ 
                          \mbox{as} \ \alpha\rightarrow+\infty.
 \end{eqnarray}
                   Combining (\ref{w1}), (\ref{5.6}) and (\ref{5.66666123}), we see
                           that (\ref{5.10}) holds.
 \par
   $Step\ 3.$ For the  fixed $\hat{\mu}\geq \mu_0$, there exist $\hat{M}\geq{{M}_0}$ and $N>0$  such that
                   ${\hat{t}}\in [T-\bar{a},T)$, $\hat{{\gamma}}_{{\hat{t}}}, \hat{{\eta}}_{{\hat{t}}}\in
                {\cal{C}}_{\hat{t}, \hat{M}}^{\hat{\mu}}$ and $|\hat{{\gamma}}_{{\hat{t}}}({\hat{t}})|\vee|\hat{{\eta}}_{{\hat{t}}}({\hat{t}})|<{\hat{M}}$  for all $\alpha\geq N$.
%
 \par
First,   noting $\varepsilon\frac{\nu T-t}{\nu
                 T}\geq\frac{\varepsilon}{1+96TL}$, 
                   by the definition of  $\Psi$,
 there exists an $\hat{M}\geq{M}_0$  that is sufficiently  large   that
           $
           \Psi(\gamma_t, \eta_t)<0
           $ for all $t\in [T-\bar{a},T]$ and $|\gamma_t(t)|\vee |\eta_t(t)|\geq{\hat{M}}$. 
           Thus, we have $|\hat{{\gamma}}_{{\hat{t}}}({\hat{t}})|\vee|\hat{{\eta}}_{{\hat{t}}}({\hat{t}})|<{\hat{M}}$.
   \par
  Next, for the fixed $\hat{M}>0$,  by (\ref{5.66666123}), we can let $N>0$ be a large number such that
$$
                         L||\hat{{\gamma}}_{{\hat{t}}}-\hat{{\eta}}_{{\hat{t}}}||_0
                         \leq
                         \frac{\tilde{m}}{4},
$$
               for all $\alpha\geq N$.
            Then we have $\hat{t}\in [T-\bar{a},T)$ for all $\alpha\geq N$. Indeed, if say $\hat{t}=T$,  we will deduce the following contradiction:
 \begin{eqnarray*}
                         \frac{\tilde{m}}{2}\leq\Psi(\hat{{\gamma}}_{{\hat{t}}},\hat{{\eta}}_{{\hat{t}}})\leq \phi(\hat{{\gamma}}_{{\hat{t}}})-\phi(\hat{{\eta}}_{{\hat{t}}})\leq
                         L||\hat{{\gamma}}_{{\hat{t}}}-\hat{{\eta}}_{{\hat{t}}}||_0\leq
                         \frac{\tilde{m}}{4}.
 \end{eqnarray*}
 \par
 $Step\ 4.$    Completion of the proof.
\par
          From above all,  for the fixed $\hat{\mu}\geq \mu_0$ in step 1 and the fixed  $\hat{M}\geq M_0$ and $N>0$ in step 3, we  find
$\hat{{\gamma}}_{{\hat{t}}}, \hat{{\eta}}_{{\hat{t}}}\in
                 {\cal{C}}_{\hat{t}, \hat{M}}^{\hat{\mu}}$   satisfying ${\hat{t}}\in [T-\bar{a},T)$ and $|\hat{{\gamma}}_{{\hat{t}}}({\hat{t}})|\vee|\hat{{\eta}}_{{\hat{t}}}({\hat{t}})|<{\hat{M}}$  for all $\alpha\geq N$
           such that
\begin{eqnarray}\label{psi4}
                   \Psi(\hat{{\gamma}}_{{\hat{t}}}, \hat{{\eta}}_{{\hat{t}}})\geq
                   \Psi(\gamma_t, \eta_t), \ (t,\gamma_t,\eta_t)\in [T-\bar{a},T]\times{\cal{C}}_{\hat{M}}^{\hat{\mu}}\times{\cal{C}}_{\hat{M}}^{\hat{\mu}}.
\end{eqnarray}
Now we consider the functional,
              for $(t,\gamma_t), (s,\eta_s)\in [\hat{t},T]\times{{\Lambda}}$,
\begin{eqnarray}\label{4.1111}
                 \Psi_{\delta}(\gamma_t,\eta_s)=W'_{1}(\gamma_t)-W'_{2}(\eta_s)-\alpha(\Upsilon^2(\gamma_{t},\hat{{\xi}}_{{\hat{t}}})+\Upsilon^2(\eta_{s},\hat{{\xi}}_{{\hat{t}}}))-\frac{1}{\delta}|s-t|^2,
\end{eqnarray}
where
%
\begin{eqnarray*}
                             {W}'_{1}(\gamma_t)&=&W_1(\gamma_t)-\varepsilon\frac{\nu T-t}{\nu
                 T}(S(\gamma_t)+|\gamma_t(t)|^2)-\varepsilon (|t-{\hat{t}}|^2
                 +||\gamma_{t}-\hat{{\gamma}}_{{\hat{t}},t}||_H^2),\\
                             {W}'_{2}(\eta_s)&=&W_2(\eta_s)+\varepsilon\frac{\nu T-s}{\nu
                 T}(S(\eta_s)+|\eta_s(s)|^2)+\varepsilon (|s-{\hat{t}}|^2
                 +||\eta_{s}-\hat{{\eta}}_{{\hat{t}}, s}||_H^2),
\end{eqnarray*}
and
$$
\hat{\xi}_{\hat{t}}=\frac{\hat{\gamma}_{\hat{t}}+\hat{\eta}_{\hat{t}}}{2}.
$$
By Lemma \ref{lemmaleft}, it has a maximum at some point $(\check{t},\check{s},\check{\gamma}_{\check{t}},\check{\eta}_{\check{s}})$ in $[\hat{t},T]\times[\hat{t},T]\times {\cal{C}}^{\hat{\mu}}_{\hat{t}, \hat{M}}\times   {\cal{C}}^{\hat{\mu}}_{\hat{t},\hat{M}}$.
By the following Lemma \ref{lemma4.4}, we have
\begin{eqnarray}\label{4.23}
\lim_{\delta\rightarrow0}\left[\frac{1}{\delta}|\check{t}-\check{s}|^2+d_\infty(\check{\gamma}_{\check{t}},\hat{\gamma}_{\hat{t}})
+d_\infty(\check{\eta}_{\check{s}},\hat{\eta}_{\hat{t}})\right]=0.
\end{eqnarray}
From  ${\hat{t}}<T$ and $|\hat{{\gamma}}_{{\hat{t}}}({\hat{t}})|\vee|\hat{{\eta}}_{{\hat{t}}}({\hat{t}})|<{\hat{M}}$  for all $\alpha\geq N$ and (\ref{4.23}),
it follows that, for every fixed $\alpha>N$,   constant $ K_\alpha>0$ exists such that
$$
            |\check{t}|\vee|\check{s}|<T, \ \    |\check{{\gamma}}_{{\check{t}}}({\check{t}})|\vee|\check{{\eta}}_{{\check{t}}}({\check{t}})|<{\hat{M}},
            \ \ \mbox{for all}    \ \ 0<\delta< K_\alpha.
$$
Now, for every $\alpha>N$ and $0<\delta< K_\alpha$, 
since $W_1$ $(\mbox{resp}., W_2)$ is a viscosity $\hat{\mu}$-subsolution (resp., $\hat{\mu}$-supsolution) to equation (\ref{hjb1}),  from  Lemma \ref{theoremsg0809827} it follows that
\begin{eqnarray}\label{vis1}
                      &&\frac{2}{\delta}(\check{t}-\check{s})
                      -\frac{\varepsilon}{\nu T}(S(\check{\gamma}_{{\check{t}}})+|\check{\gamma}_{{\check{t}}}({\check{t}})|^2)
                      +\varepsilon|\check{\gamma}_{{\check{t}}}({\check{t}})
                      -\hat{\gamma}_{{\hat{t}}}({\hat{t}})|^2
                      +2\varepsilon({\check{t}}-{\hat{t}})
                      \nonumber\\
                      &&
                      +{\mathbf{H}}\left(\check{\gamma}_{{\check{t}}},
                      \varepsilon\frac{\nu T-{\check{t}}}{\nu T}(\partial_xS(\check{\gamma}_{{\check{t}}})
                       +2 \check{\gamma}_{{\check{t}}}({\check{t}}))
                      +4\alpha(\check{\gamma}_{{\check{t}}}({\check{t}})-\hat{\xi}_{\hat{t}}(\hat{t}))
                      +\alpha\partial_xS^{\hat{\xi}_{\hat{t}}}(\check{\gamma}_{{\check{t}}})
                                    \right)
                                     \geq c;
\end{eqnarray}
and
\begin{eqnarray}\label{vis2}
                      &&\frac{2}{\delta}(\check{t}-\check{s})
                     +\frac{\varepsilon}{\nu T}(S(\check{\eta}_{{\check{s}}})+|\check{\eta}_{{\check{s}}}({\check{s}})|^2)
                      -\varepsilon|\check{\eta}_{{\check{s}}}({\check{s}})-\hat{\eta}_{{\hat{t}}}({\hat{t}})|^2-2\varepsilon({\check{s}}-{\hat{t}})
                      \nonumber\\
                      &&+{\mathbf{H}}\left(\check{\eta}_{{\check{s}}},  -\varepsilon\frac{\nu T-{\check{s}}}{\nu T}
                     (\partial_xS(\check{\eta}_{{\check{t}}})
                       +2 \check{\eta}_{{\check{s}}}(\check{{s}}))
                       -4\alpha(\check{\eta}_{{\check{s}}}({\check{s}})-\hat{\xi}_{\hat{t}}(\hat{t}))
                      -\alpha\partial_xS^{\hat{\xi}_{\hat{t}}}(\check{\eta}_{{\check{s}}})
                                     \right)
                                     \leq 0.
\end{eqnarray}
 Combining (\ref{vis1}) and (\ref{vis2}), and letting $\delta\rightarrow0$,   we obtain
 \begin{eqnarray}\label{vis112}
                     &&c+ \frac{\varepsilon}{\nu T}(S(\hat{{\gamma}}_{{\hat{t}}})+|\hat{{\gamma}}_{{\hat{t}}}({\hat{t}})|^2+S(\hat{{\eta}}_{{\hat{t}}})+|\hat{{\eta}}_{{\hat{t}}}({\hat{t}})|^2)\nonumber\\
                     &\leq&
                                     {\mathbf{H}}\left(\hat{{\gamma}}_{{\hat{t}}},  2\alpha(\hat{{\gamma}}_{{\hat{t}}}({\hat{t}})-\hat{{\eta}}_{{\hat{t}}}({\hat{t}})) +\varepsilon\frac{\nu T-{\hat{t}}}{\nu T}(\partial_xS(\hat{{\gamma}}_{{\hat{t}}})
                                     +2\hat{{\gamma}}_{{\hat{t}}}({\hat{t}}))
                                      +\alpha\partial_xS^{\hat{\xi}_{\hat{t}}}(\hat{\gamma}_{{\hat{t}}})
                                     \right)\nonumber\\
                                     &&-{\mathbf{H}}\left(\hat{{\eta}}_{{\hat{t}}},
                                     2\alpha(\hat{{\gamma}}_{{\hat{t}}}({\hat{t}})-\hat{{\eta}}_{{\hat{t}}}({\hat{t}})) -\varepsilon\frac{\nu T-{\hat{t}}}{\nu T}(\partial_xS(\hat{{\eta}}_{{\hat{t}}})
                                     +2\hat{{\eta}}_{{\hat{t}}}({\hat{t}}))
                                      -\alpha\partial_xS^{\hat{\xi}_{\hat{t}}}(\hat{\eta}_{{\hat{t}}})\right).
\end{eqnarray}
  On the other hand, by  a simple calculation we obtain
\begin{eqnarray}\label{v4}
                &&                               {\mathbf{H}}\left(\hat{{\gamma}}_{{\hat{t}}},  2\alpha(\hat{{\gamma}}_{{\hat{t}}}({\hat{t}})-\hat{{\eta}}_{{\hat{t}}}({\hat{t}})) +\varepsilon\frac{\nu T-{\hat{t}}}{\nu T}(\partial_xS(\hat{{\gamma}}_{{\hat{t}}})
                                     +2\hat{{\gamma}}_{{\hat{t}}}({\hat{t}}))
                                      +\alpha\partial_xS^{\hat{\xi}_{\hat{t}}}(\hat{\gamma}_{{\hat{t}}})
                                     \right)\nonumber\\
                                     &&-{\mathbf{H}}\left(\hat{{\eta}}_{{\hat{t}}},
                                     2\alpha(\hat{{\gamma}}_{{\hat{t}}}({\hat{t}})-\hat{{\eta}}_{{\hat{t}}}({\hat{t}}))
                 -\varepsilon\frac{\nu T-{\hat{t}}}{\nu T}(\partial_xS(\hat{{\eta}}_{{\hat{t}}})
                                     +2\hat{{\eta}}_{{\hat{t}}}({\hat{t}}))
                                      -\alpha\partial_xS^{\hat{\xi}_{\hat{t}}}(\hat{\eta}_{{\hat{t}}})\right)
\nonumber\\
                &\leq&\sup_{u\in U}(J_{1}+J_{2}),
\end{eqnarray}
            where
\begin{eqnarray}\label{j1}
                               J_{1}
                               &=&\left( {F}(\hat{{\gamma}}_{{\hat{t}}},u),2\alpha(\hat{{\gamma}}_{{\hat{t}}}({\hat{t}})-\hat{{\eta}}_{{\hat{t}}}({\hat{t}}))
                               +\varepsilon\frac{\nu T-{\hat{t}}}{\nu T}(\partial_xS(\hat{{\gamma}}_{{\hat{t}}})
                                     +2\hat{{\gamma}}_{{\hat{t}}}({\hat{t}}))
                                     +\alpha\partial_xS^{\hat{\xi}_{\hat{t}}}(\hat{\gamma}_{{\hat{t}}})\right)_{R^d}\nonumber\\
                                            &&  -\left( {F}(\hat{{\eta}}_{{\hat{t}}},u),2\alpha(\hat{{\gamma}}_{{\hat{t}}}({\hat{t}})-\hat{{\eta}}_{{\hat{t}}}({\hat{t}}))
                                            -\varepsilon\frac{\nu T-{\hat{t}}}{\nu T}(\partial_xS(\hat{{\eta}}_{{\hat{t}}})
                                     +2\hat{{\eta}}_{{\hat{t}}}({\hat{t}}))
                                     -\alpha\partial_xS^{\hat{\xi}_{\hat{t}}}(\hat{\eta}_{{\hat{t}}})\right)_{R^d}\nonumber\\
                                  &\leq&4\alpha{L}|\hat{{\gamma}}_{{\hat{t}}}({\hat{t}})-\hat{{\eta}}_{{\hat{t}}}({\hat{t}})|
                                 \times||\hat{{\gamma}}_{{\hat{t}}}-\hat{{\eta}}_{{\hat{t}}}||_0
                                 +6\varepsilon \frac{\nu T-{\hat{t}}}{\nu T} L|\hat{{\gamma}}_{{\hat{t}}}({\hat{t}})|(1+||\hat{{\gamma}}_{{\hat{t}}}||_0
                                             )\nonumber\\
                                            && +6\varepsilon \frac{\nu T-{\hat{t}}}{\nu T} L|\hat{{\eta}}_{{\hat{t}}}({\hat{t}})|(1+||\hat{{\eta}}_{{\hat{t}}}||_0
                                             );
\end{eqnarray}
and
\begin{eqnarray}\label{j3}
                                 J_{2}=q(\hat{{\gamma}}_{{\hat{t}}},u)-
                                 q(\hat{{\eta}}_{{\hat{t}}},u)\leq
                                 L||\hat{{\gamma}}_{{\hat{t}}}-\hat{{\eta}}_{{\hat{t}}}||_0.
\end{eqnarray}
 Combining (\ref{vis112})-(\ref{j3}), we obtain
 \begin{eqnarray}\label{vis122}
                    c
                                             &\leq&
                      -\frac{\varepsilon}{\nu T}(S(\hat{{\gamma}}_{{\hat{t}}})+|\hat{{\gamma}}_{{\hat{t}}}({\hat{t}})|^2
                     +S(\hat{{\eta}}_{{\hat{t}}})+|\hat{{\eta}}_{{\hat{t}}}({\hat{t}})|^2)
                      +2\alpha{L}(|\hat{{\gamma}}_{{\hat{t}}}({\hat{t}})-\hat{{\eta}}_{{\hat{t}}}({\hat{t}})|^2+||\hat{{\gamma}}_{{\hat{t}}}-\hat{{\eta}}_{{\hat{t}}}||_0^2)\nonumber\\
                                 &&
                                 +L||\hat{{\gamma}}_{{\hat{t}}}-\hat{{\eta}}_{{\hat{t}}}||_0
                               +12\varepsilon \frac{\nu T-{\hat{t}}}{\nu T}L(1+||\hat{{\gamma}}_{{\hat{t}}}||_0^2
                                             +||\hat{{\eta}}_{{\hat{t}}}||_0^2
).
\end{eqnarray}
Recalling $
         \nu=1+\frac{1}{96TL}
$ and $\bar{a}=\frac{1}{96L}\wedge\frac{T}{2}$, by (\ref{s0}),  we have
\begin{eqnarray*}
                     c
                                 \leq
                             2\alpha L(|\hat{{\gamma}}_{{\hat{t}}}({\hat{t}})-\hat{{\eta}}_{{\hat{t}}}({\hat{t}})|^2
                              +||\hat{{\gamma}}_{{\hat{t}}}-\hat{{\eta}}_{{\hat{t}}}||_0^2)
                              +L||\hat{{\gamma}}_{{\hat{t}}}-\hat{{\eta}}_{{\hat{t}}}||_0
                              +\frac{\varepsilon}{\nu
                              T}.
\end{eqnarray*}
Then,  letting $\alpha\rightarrow\infty$, by (\ref{s0}), (\ref{5.3}) and (\ref{5.10}), the following contradiction is induced:
\begin{eqnarray*}
                    c
                      \leq\frac{c}{2}.
\end{eqnarray*}
 The proof is now complete.
 \ \ $\Box$

 To complete the previous proof, it remains to state and prove the following lemma.
 \begin{lemma}\label{lemma4.4}\ \ The maximum point $(\check{t},\check{s},\check{\gamma}_{\check{t}},\check{\eta}_{\check{s}})$
of $\Psi_{\delta}(\gamma_t,\eta_s)
        $ defined by (\ref{4.1111}) in  $[\hat{t},T]\times [\hat{t},T]\times {\cal{C}}^{\hat{\mu}}_{\hat{t}, \hat{M}}\times   {\cal{C}}^{\hat{\mu}}_{\hat{t}, \hat{M}}$ satisfies   condition (\ref{4.23}).
\end{lemma}
\par
   {\bf  Proof  }. \ \
  By the definition of the maximum point $(\check{t},\check{s},\check{\gamma}_{\check{t}},\check{\eta}_{\check{s}})$, we have
 \begin{eqnarray*}
                        2\Psi_{\delta}(\check{\gamma}_{\check{t}},\check{\eta}_{\check{s}})
            \geq  \Psi_{\delta}(\check{\gamma}_{\check{t}},\check{\gamma}_{\check{t}})+\Psi_{\delta}(\check{\eta}_{\check{s}},\check{\eta}_{\check{s}}).
 \end{eqnarray*}
 This implies that
 \begin{eqnarray*}
                         \frac{2}{\delta}|\check{t}-\check{s}|^2
                         \leq|W'_1(\check{\gamma}_{\check{t}})-W'_1(\check{\eta}_{\check{s}})|
                                   +|W'_2(\check{\gamma}_{\check{t}})-W'_2(\check{\eta}_{\check{s}})|.
 \end{eqnarray*}
 Letting $\delta\rightarrow0$, we have
 $$
 |\check{t}-\check{s}|\rightarrow0 \ \ \mbox{as}\  \ \delta\rightarrow0.
 $$
  Without loss of generality, we assume $\check{t}\leq \check{s}$.
  By the definition  of $\Psi_{\delta}$ and (\ref{jias5}), we get that
\begin{eqnarray*}
                 &&\Psi_{\delta}(\check{\gamma}_{\check{t}},\check{\eta}_{\check{s}})
                  \geq\Psi_{\delta}(\hat{\gamma}_{\hat{t}},\hat{\eta}_{\hat{t}})=\Psi(\hat{\gamma}_{\hat{t}},\hat{\eta}_{\hat{t}})
                 \geq \Psi(\check{\gamma}_{\check{t},\check{s}},\check{\eta}_{\check{s}})\\
                 &=& W_1(\check{\gamma}_{\check{t},\check{s}})- {W}_2(\check{\eta}_{\check{s}})-\frac{\alpha}{2}\Upsilon^2(\check{\gamma}_{\check{t},\check{s}},\check{\eta}_{\check{s}})
                 -\varepsilon\frac{\nu T-\check{s}}{\nu T}(S(\check{\gamma}_{\check{t},\check{s}})+S(\check{\eta}_{\check{s}})+|\check{\gamma}_{\check{t}}(\check{t})|^2+|\check{\eta}_{\check{s}}(\check{s})|^2)
                 \\
                 &\geq&\Psi_{\delta}(\check{\gamma}_{\check{t}},\check{\eta}_{\check{s}})+\frac{1}{\delta}|\check{t}-\check{s}|^2-L(1+\hat{M})|\check{t}-\check{s}|
                                     +\varepsilon (|\check{t}-{\hat{t}}|^2+|\check{s}-{\hat{t}}|^2
                 +||\check{\gamma}_{\check{t}}-\hat{{\gamma}}_{{\hat{t}},\check{t}}||_H^2
                 +||\check{\eta}_{\check{t}}-\hat{{\eta}}_{{\hat{t}},\check{s}}||_H^2).
\end{eqnarray*}
Letting $\delta\rightarrow0$, we obtain that
$$
\frac{1}{\delta}|\check{t}-\check{s}|^2
                                     +\varepsilon (|\check{t}-{\hat{t}}|^2+|\check{s}-{\hat{t}}|^2
                 +||\check{\gamma}_{\check{t}}-\hat{{\gamma}}_{{\hat{t}},\check{t}}||_H^2
                 +||\check{\eta}_{\check{s}}-\hat{{\eta}}_{{\hat{t}},\check{s}}||_H^2)\rightarrow0\ \ \mbox{as}\  \ \delta\rightarrow0.
$$
 On the other hand, since $\check{\gamma}_{\check{t}},\check{\eta}_{\check{s}}\in {\cal{C}}^{\hat{\mu}}_{\hat{t}, \hat{M}}$ and $|\check{t}-{\hat{t}}|+|\check{s}-{\hat{t}}|\rightarrow0$ as $\delta\rightarrow0$, we may assume $d_\infty(\check{\gamma}_{\check{t}},\bar{\gamma}_{{\hat{t}}})
 +d_\infty(\check{\eta}_{\check{s}},\bar{\eta}_{{\hat{t}}})\rightarrow0$
                 for some $\bar{\gamma}_{{\hat{t}}},\bar{\eta}_{{\hat{t}}}\in {\cal{C}}^{\hat{\mu}}_{\hat{t}, \hat{M}}$. Then we have that
\begin{eqnarray*}
                   ||\check{\gamma}_{\check{t}}-\bar{{\gamma}}_{{\hat{t}},\check{t}}||_H^2
                 +||\check{\eta}_{\check{s}}-\bar{{\eta}}_{{\hat{t}},\check{s}}||_H^2\rightarrow0\ \ \mbox{as}\  \ \delta\rightarrow0.
\end{eqnarray*}
Therefore, $\hat{\gamma}_{{\hat{t}}}=\bar{\gamma}_{{\hat{t}}},\hat{\eta}_{{\hat{t}}}=\bar{\eta}_{{\hat{t}}}$,
and we get that (\ref{4.23}) holds true.  The proof is now complete. \ \ $\Box$

\section{Appendix}  \label{RDS}

\par
In this Appendix, we prove $(\hat{\Lambda}^t, d_{\infty})$ is a complete metric space.
\begin{lemma}\label{lemma2.1111}
 $(\hat{\Lambda}^t, d_{\infty})$ is a complete metric space for every $t\in [0,T)$. 
\end{lemma}
\par
{\bf  Proof}. \ \
Assume $\{\gamma^n_{t_n}\}_{n\geq0}$ is a cauchy sequence in $(\hat{\Lambda}^t, d_{\infty})$, then for any $\varepsilon>0$,  there exists $N(\varepsilon)>0$ such that, for all $m,n\geq N(\varepsilon)$, we have
\begin{eqnarray*}\label{jialemma111000}
d_\infty(\gamma^n_{t_n},\gamma^m_{t_m})=|t_n-t_m|+\sup_{0\leq s\leq T}|\gamma^n_{t_n}(s\wedge t_n)-\gamma^m_{t_m}(s\wedge t_m)|<\varepsilon.
\end{eqnarray*}
Therefore, there exists  $\hat{t}\in [t,T]$ such that
$\lim_{n\rightarrow\infty}t_n=\hat{t}.$
 Moreover, for all $s\in [0,T]$,
\begin{eqnarray}\label{jialemma}
    |\gamma^n_{t_n}(s\wedge t_n)-\gamma^m_{t_m}(s\wedge t_m)|<\varepsilon, \  (\forall m,n\geq N(\varepsilon)).
\end{eqnarray}
For fixed $s\in [0,T]$, we see that $\{\gamma^n_{t_n}(t_n\wedge s)\}$ is a cauchy sequence, thereby the limit $\lim_{n\rightarrow\infty}\gamma^n_{t_n}(t_n\wedge s)$  exists and denoted by $\gamma_{T}(s)$. Letting $m\rightarrow\infty$ in (\ref{jialemma}), we obtain that
\begin{eqnarray*}
    |\gamma_{T}(s)-\gamma^n_{t_n}(s\wedge t_n)|\leq\varepsilon, \  (\forall  \ n\geq N(\varepsilon)).
\end{eqnarray*}
 Taking the  supremum over $s\in
             [0,T]$, we get
             \begin{eqnarray}\label{jialemma3}
    \sup_{s\in[0,T]}|\gamma_{T}(s)-\gamma^n_{t_n}(s\wedge t_n)|\leq\varepsilon, \  (\forall \ n\geq N(\varepsilon)).
\end{eqnarray}
We claim that $\gamma_{T}(s)=\gamma_T(\hat{t})$  for all $s\in(\hat{t},T]$. In fact, if there exists a subsequence $\{t_{n_l}\}_{l\geq0}$ of $\{t_{n}\}_{n\geq0}$ such that $\{t_{n_l}\}_{l\geq0}\leq\hat{t}$, then we have, for every
 $s\in(\hat{t},T]$,
 $$
         \gamma_T(s)=\lim_{n\rightarrow\infty}\gamma^n_{t_n}(s\wedge t_n)=\lim_{n\rightarrow\infty}\gamma^n_{t_n}(t_n)=\lim_{l\rightarrow\infty}\gamma^n_{t_{n_l}}(t_{n_l})=\lim_{l\rightarrow\infty}\gamma^n_{t_{n_l}}(t_{n_l}\wedge \hat{t})=\lim_{n\rightarrow\infty}\gamma^n_{t_{n}}(t_{n}\wedge \hat{t})=\gamma_T(\hat{t}).
 $$
 Otherwise, we may assume $\{t_{n}\}_{n\geq0}>\hat{t}$. Letting  $s=t_m$ and $m\rightarrow\infty$ in (\ref{jialemma}), we obtain, for all $s\in (\hat{t},T]$,
 $$
  |\gamma^n_{t_n}(\hat{t})-\gamma_T(s)|\leq\varepsilon, \  (\forall n\geq N(\varepsilon)).
 $$
 Letting $n\rightarrow\infty$, we have
 $$
 |\gamma_{T}(\hat{t})-\gamma_T(s)|\leq\varepsilon,\ \mbox{for all}\ \ s\in(\hat{t},T].
 $$
Then, by (\ref{jialemma3}) we obtain
\begin{eqnarray}\label{jialemma23}
d_\infty(\eta_{\hat{t}},\gamma^n_{t_n})\rightarrow0 \ \mbox{as} \  n\rightarrow\infty.
\end{eqnarray}
Here we let $ \eta_{\hat{t}}$ denote $\gamma_T|_{[0,\hat{t}]}$.
Now we prove $\eta_{\hat{t}}\in \hat{\Lambda}^t$. First, we prove $\eta_{\hat{t}}$ is right-continuous.
For every $0\leq s<\hat{t}$ and $0<\delta\leq \hat{t}-s$, we have
$$
|\eta_{\hat{t}}(s+\delta)-\eta_{\hat{t}}(s)|\leq |\gamma_{T}(s+\delta)-\gamma^n_{t_n}((s+\delta)\wedge t_n)|+|\gamma^n_{t_n}((s+\delta)\wedge t_n)-\gamma^n_{t_n}(s\wedge t_n)|+|\gamma^n_{t_n}(s\wedge t_n)-\gamma_{T}(s)|.
$$
For every $\varepsilon>0$, by (\ref{jialemma3}), there exists $n>0$ independent of $\delta$, which is large enough such that
$$
|\gamma_{T}(s+\delta)-\gamma^n_{t_n}((s+\delta)\wedge t_n)|+|\gamma^n_{t_n}(s\wedge t_n)-\gamma_{T}(s)|< \frac{\varepsilon}{2}.
$$
For the fixed $n$, since $\gamma^n_{t_n}\in \hat{\Lambda}^t$, there exists a constant $0<\Delta\leq \hat{t}-s$ such that, for all $0\leq\delta<\Delta$,
$$
  |\gamma^n_{t_n}((s+\delta)\wedge t_n)-\gamma^n_{t_n}(s\wedge t_n)|< \frac{\varepsilon}{2}.
$$
Then $
|\eta_{\hat{t}}(s+\delta)-\eta_{\hat{t}}(s)|< \varepsilon
$ for all $0\leq\delta<\Delta$. Next, let us prove  $\eta_{\hat{t}}$  has left limit in $(0,\hat{t}]$. For every $0< s\leq\hat{t}$ and $0\leq s_1,s_2<s$, we have
$$
|\eta_{\hat{t}}(s_1)-\eta_{\hat{t}}(s_2)|\leq |\gamma_{T}(s_1)-\gamma^n_{t_n}(s_1\wedge t_n)|+|\gamma_{T}(s_2)-\gamma^n_{t_n}(s_2\wedge t_n)|+|\gamma^n_{t_n}(s_1\wedge t_n)-\gamma^n_{t_n}(s_2\wedge t_n)|.
$$
For every $\varepsilon>0$, by (\ref{jialemma3}), there exists $n>0$ be large enough such that
$$
|\gamma_{T}(s_1)-\gamma^n_{t_n}(s_1\wedge t_n)|+|\gamma_{T}(s_2)-\gamma^n_{t_n}(s_2\wedge t_n)|< \frac{\varepsilon}{2}.
$$
For the fixed $n$, if $t_n<s$, we can let $\Delta>0$ be small enough such that $t_n<s-\Delta$, then for all  $s_1,s_2\in [s-\Delta,s)$,
$$
|\gamma^n_{t_n}(s_1\wedge t_n)-\gamma^n_{t_n}(s_2\wedge t_n)|=|\gamma^n_{t_n}(t_n)-\gamma^n_{t_n}(t_n)|=0;
$$
if $t_n\geq s$, since $\gamma^n_{t_n}\in \hat{\Lambda}^t$, there exists a constant $\Delta>0$ such that, for all $s_1,s_2\in [s-\Delta,s)$,
$$
|\gamma^n_{t_n}(s_1\wedge t_n)-\gamma^n_{t_n}(s_2\wedge t_n)|=|\gamma^n_{t_n}(s_1)-\gamma^n_{t_n}(s_2)|<\frac{\varepsilon}{2}.
$$
Then there exists a constant $\Delta>0$ such that $
|\eta_{\hat{t}}(s_1)-\eta_{\hat{t}}(s_2)|< \varepsilon
$ for all  $s_1,s_2\in [s-\Delta,s)$. The proof is now complete. \ \ $\Box$

\par

\end{document}